\documentclass[11pt]{article}
\def\WHO{nobody}
\def\version{16.10.2019 - REVISION 17.4.2020}\def\users{}  %
\def\users{final-layout}  
\usepackage{amsmath,amsfonts,amssymb,color}
\usepackage{mathrsfs} 
\usepackage{epsfig}
\usepackage{psfrag}
\usepackage{upgreek} 
\usepackage{cite}
\usepackage{url}
\usepackage{import}
 \usepackage[symbol]{footmisc}

\usepackage{setspace}
\usepackage{fullpage}

\newtheorem{theorem}{Theorem}[section]
\newtheorem{definition}[theorem]{Definition}
\newtheorem{proposition}[theorem]{Proposition}
\newtheorem{remark}[theorem]{Remark}

\definecolor{darkgreen}{rgb}{0.2,0.7,0.2}

\numberwithin{equation}{section}
\usepackage{ifthen}
\ifthenelse{\equal{\users}{final-layout}}{}{
\usepackage{fancyhdr}
\pagestyle{fancy}
\headheight=28pt\headwidth=16.5cm
\definecolor{gray}{gray}{0.5}
\rhead{\color{gray}Explicit elastodynamics\\
T.Roub\'\i\v cek and C. Tsogka}
\chead{}
\lhead{Version \version, file: \jobname.tex, \underline{\Large\color{red}\WHO's working on} \\
compiled:
\number\day.\number\month.\number\year\ at
\the\hour:\ifnum\minute<10 0\fi\the\minute\ CET }
}

\definecolor{labelkey}{rgb}{1.,.2,0.}

\newcount\hour \newcount\minute
\hour=\time
\divide \hour by 60
\minute=\time
\loop \ifnum \minute > 59 \advance \minute by -60 \repeat

\usepackage[normalem]{ulem}
\definecolor{brown}{rgb}{0.5,0,0}
\ifthenelse{\equal{\users}{final-layout}}{
    \newcommand{\REPLACE}[2]{#2}
    \newcommand{\INSERT}[1]{#1}
    \newcommand{\DELETE}[1]{}
    \newcommand{\CHECK}[1]{#1}
    \newcommand{\COMMENT}[1]{}
    \newcommand{\COLOR}[1]{#1}
    \newcommand{\mman}[1]{#1}
     \newcommand{\cman}[1]{#1}
    \newcommand{\TINY}[1]{}
    \newcommand{\MARGINOTE}[1]{}
}
{\newcommand{\REPLACE}[2]{{\color{brown}\sout{#1}\uline{#2}\color{black}}}
\newcommand{\INSERT}[1]{{\color{blue}\uuline{#1}\color{black}}}
 \newcommand{\DELETE}[1]{{\color{brown}\sout{#1}\color{black}}}
 \newcommand{\CHECK}[1]{\color{brown}\uwave{#1}\color{black}}
 \newcommand{\COMMENT}[1]{{\color{red}\uuline{#1}\color{black}}}
 \newcommand{\COLOR}[2][black]{{\color{#1}{#2}}}  
 \newcommand{\mman}[1]{{\color{blue}#1}}
\newcommand{\cman}[1]{{\color{darkgreen}#1}}

 \newcommand{\TINY}[1]{{\tiny#1}}
 \newcommand{\MARGINOTE}[1]{\marginpar{\color{red}\tiny\texttt{#1}}}
}

\renewcommand\dot[1]{\mathchoice
                 {{\buildrel{\hspace*{.1em}\text{\LARGE.}}\over{#1}}}
                 {{\buildrel{\hspace*{.1em}\text{\Large.}}\over{#1}}}
                 {{\buildrel{\hspace*{.1em}\text{\large.}}\over{#1}}}
                 {{\buildrel{\hspace*{.1em}\text{\large.}}\over{#1}}}}
\newcommand\DT{\dot}
\newcommand\DDT[1]{\mathchoice
   {{\buildrel{\hspace*{.1em}\text{\LARGE.\hspace*{-.1em}.}}\over{#1}}}
   {{\buildrel{\hspace*{.1em}\text{\Large.\hspace*{-.1em}.}}\over{#1}}}
   {{\buildrel{\hspace*{.1em}\text{\large.\hspace*{-.1em}.}}\over{#1}}}
   {{\buildrel{\hspace*{.1em}\text{\large.\hspace*{-.1em}.}}\over{#1}}}}

\newcommand{\wt}[1]{\mathchoice
     {\text{\small$\widetilde{\text{\normalsize$#1$}}\hspace*{.03em}$}}
                    {\text{\small$\widetilde{\text{\normalsize$#1$}}$}}
                    {\widetilde{#1\hspace*{-.02em}}\hspace*{.03em}}
                    {\widetilde{#1}}}

\def\R{{\mathbb R}}
\newcommand\calE{\mathcal E}
\newcommand\calH{\mathcal H}
\newcommand\calS{\mathcal S}
\newcommand\calU{\mathcal U}
\newcommand\calZ{\mathcal Z}
\newcommand\calX{\mathcal X}
\newcommand\bbC{\mathbb C}
\newcommand\bbD{\mathbb D}

\newcommand\bbB{\mathbb B}

\newcommand{\lineunder}[2]{\LU{\begin{array}[t]{c}\underbrace{#1}\vspace*{.5em}\end{array}}{\mbox{\footnotesize\rm #2}}}
\newcommand{\linesunder}[3]{\LSU{\begin{array}[t]{c}\underbrace{#1}\vspace*{.5em}\end{array}}{\mbox{\footnotesize\rm #2}}{\mbox{\footnotesize\rm #3}}}

\newcommand{\LU}[2]{\begin{array}[t]{c}#1\vspace*{-1em}\\_{#2}\end{array}}
\newcommand{\LSU}[3]{\begin{array}[t]{c}#1\vspace*{-1em}\\_{#2}\vspace*{-.5em}\\_{#3}\end{array}}

\renewcommand\d{\mathrm d}

\newcommand{\KK}{k+1}
\newcommand{\KKK}{k}
\newcommand{\TT}{{\mathcal T}}

\newcommand{\barr}{\begin{array}}
\newcommand{\earr}{\end{array}}

\newcommand{\ccc}{\color{black}}
\newcommand{\eee}{\color{black}}

\allowdisplaybreaks

\def\In{{\in}}
\newcommand{\eq}[1]{\eqref{#1}}

\begin{document}
\title{Staggered
explicit-implicit time-discretization
for elastodynamics with dissipative internal variables}
\author{Tom\'{a}\v{s} Roub\'\i\v{c}ek\footnotemark[1]  \footnotemark[2] and Chrysoula Tsogka\footnotemark[3] \footnotemark[4]}

\footnotetext[1]{Charles University,
Mathematical Institute, Sokolovsk\'a 83, CZ-186~75~Praha~8,
Czech Republic}
\footnotetext[2]{Institute of
Thermomechanics of the Czech Acad. Sci., Dolej\v skova~5,
CZ--182~08 Praha 8, Czech Rep.}
\footnotetext[3]{Department of Applied Mathematics,
  University of California, Merced, 5200 North Lake road, Merced, CA 95343, USA}
\footnotetext[4]{Institute of Applied and Computational Mathematics,
  Foundation for Research and Technology - Hellas, Nikolaou Plastira 100, Vassilika Vouton, GR-700 13
Heraklion, Crete, Greece}
  
\maketitle \date{}
\begin{abstract}
 An extension of the two-step staggered time discretization of
    linear elastodynamics in stress-velocity form
    to systems involving internal variables 
subjected to a possibly non-linear dissipative evolution
 is proposed. The original scheme is thus
    enhanced by another step for the internal variables which, in general,
    is implicit, although even this step might be explicit if no
    spatial gradients of the internal variables are involved.
    Using an abstract Banach-space formulation,
a-priori estimates and convergence
 are  proved under a CFL condition. The developed
  three-step scheme
  finds applications in various
  problems of continuum mechanics at small strain.
  Here, we consider in particular 
  plasticity, viscoelasticity (creep), diffusion in poroelastic media,
    and damage.
\end{abstract}



{\bf Keywords} {Elastodynamics, explicit
  discretization, \mman{fractional steps}, 
mixed finite-element method, plasticity, \mman{creep}, poro-elasticity,
\mman{damage}.}

\section{Introduction \DELETE{-- mere linear elastodynamics}}
\cman{In computational mechanics one can distinguish
two main classes of time-dependent problems\mman{, quasistatic and dynamic.
  Focusing on the latter one, one can further distinguish \cman{two} other
  cases: \cman{(i)} Low-frequency regimes, 
which are typically related with vibrations of structures and where
the energy is not dominantly transmitted through space. \cman{In this case} implicit
time-discretization is relatively efficient, even though  large
system\cman{s} of algebraic equations are to be solved at each time step.}
\cman{(ii)} 
High-frequency
\mman{regimes} which
\mman{arise} typically \mman{within} wave
propagation
\mman{and for which only explicit time-discretization\cman{s}  are
  reasonably efficient\cman{,} in particular in three-dimensions. 
  These explicit methods can essentially be used in hyperbolic problems,
  as mere elastodynamics (treated here) or the elasto-magnetic Maxwell system
  or some conservation laws.
Yet, many applications need to combine the convervative hyperbolic problems
with various dissipative processes of a parabolic character, involving
typically some internal variables. For parabolic problems, however,
explicit methods are known to be problematic due to severe time-step
restrictions. Therefore, we propose and analyse} 
an explicit/implicit scheme\mman{, using the fractional-step (also called
  staggered) technique}.
In fact, the proposed 
method becomes completely explicit in the absence of gradients of internal
variables, as} \mman{it can be in plasticicity or viscoelasticity
  as in Sect.~\ref{sect-plast}, or an interfacial variant of damage
  models (so-called delamination) as in \cite{RKVPZ??DACM}.}

Our starting point is the linear elastodynamic problem: Find the displacement
$u:[0,T]\times\varOmega\to\R^d$ satisfying\hspace*{-.1em}
%
\begin{subequations}\label{IBVP}\begin{align}\label{IBVP1}
    &&&  \varrho\DDT u-{\rm div}\,\bbC e(u)=f&&\text{on }\ \varOmega\
    \text{ for }\ t\in\mman{(}0,T],&&
  \\&&&\label{IBVP2}
    [\bbC e(u)]\vec{n}+\bbB u=g&&\text{on }\ \varGamma\
    \text{ for }\ t\in\mman{(}0,T],
  \\&&&\label{IBVP3}
  u|_{t=0}=u_0,\ \ \ \DT u|_{t=0}=v_0&&\text{on }\ \varOmega,
\end{align}\end{subequations}
\cman{for $T>0$ a fixed time horizon. Here $\varOmega\subset\R^d$ is a bounded Lipschitz domain, $d=2$ or $3$, 
$\varGamma$ is its boundary, and $\vec{n}$ the unit outward normal. The dot-notation stands for the time derivative.
In  (\ref{IBVP}), $\varrho>0$ denotes the mass density, $\bbC$ is the 
elasticity tensor which is symmetric and positive definite, $e(u)$ denotes the \mman{small-}strain tensor defined as
$e(u)=\frac12(\nabla u)^\top{+}\frac12\nabla u$, and 
$\bbB$ is a symmetric positive semidefinite 2nd-order tensor determining the elastic support on the boundary. The terms appearing in the second member of (\ref{IBVP}) are,  
the bulk force $f$ and the surface loading $g$. In (\ref{IBVP3}),  
$u_0$ denotes the 
initial displacement and $v_0$ the initial velocity.}
\cman{For a more compact notation, we write the initial-boundary-value
problem \eqref{IBVP} in the following abstract form }
\begin{align}\label{IVP}
  \mathcal T'\DDT u+\mathcal W' u=\mathcal F_u'(t)\ \ 
\text{ for }\ t\in\mman{(}0,T],\ \ u|_{t=0}=u_0,\ \ \ \DT u|_{t=0}=v_0.
\end{align}
Here $\mathcal T$ is the kinetic energy, $\mathcal W$ is the stored
energy, and $\mathcal F$ is the external force, while
$(\cdot)'$ denotes the G\^ateaux derivative. In the context
of \eqref{IBVP}, we have 
$$
\mathcal T(\ccc\DT u\eee)=\int_\varOmega\frac12\varrho| \ccc \DT u\eee|^2\,\d x,\ \quad 
\mathcal W(u)=\int_\varOmega\frac12\bbC e(u){:}e(u)\,\d x
+\int_\varGamma\frac12\bbB u\cdot u\,\d S
$$
and
$$
\mathcal F(t,u)=\int_\varOmega f(t){\cdot}u\,\d x
+\int_\varGamma g(t){\cdot}u \,\d S.
$$
Thus $\mathcal F_u'(t)$ is a linear functional, let us denote it shortly by
$F(t)$.

\cman{For high-frequency wave propagation problems, implicit time
  discretizations are computationally expensive, especially in the three
  dimensional setting. Therefore, we focus our attention on explicit methods. } 
The simplest explicit scheme is the following second-order finite 
difference scheme
\begin{align}\label{elast-dyn-explicit}
\mathcal T'\frac{u_{\tau h}^{k+1}-2u_{\tau h}^{k}+u_{\tau h}^{k-1}}{\tau^2}
+\mathcal W_h' u_{\tau h}^{k}=F_h(k\tau).
\end{align}
Here $\tau>0$ denotes the time step, and 
$\mathcal W_h$ (resp.\ $F_h$) denote some discrete approximations of the
respective continuous functionals
obtained by a suitable finite-element method (FEM) with mesh size $h>0$.
\mman{For simplicity, we assume the mass density  constant
  (or at least piecewise constant in space)  so that the kinetic energy
  $\mathcal T$ does not need any numerical
  approximation.}
In particular, a numerical approximation 
leading to a diagonal the mass matrix $\mathcal T'$ in \eq{elast-dyn-explicit}, typically referred to as mass lumping,  
is an important ingredient so as to obtain 
efficient explicit methods.
\ccc Here we will consider that $u$ is discretized in an element-wise constant way
so that $\mathcal T'$ leads to a diagonal matrix even without any
approximation. 
\eee
Multiplying \eqref{elast-dyn-explicit} by $\frac{u_{\tau h}^{k+1}-u_{\tau h}^{k}}{\tau}$ and using the scheme, it is easy to show that the 
energy preserved is  
$$
\frac12\langle
\ccc\TT'\eee
\frac{u_{\tau h}^{k+1}-u_{\tau h}^{k}}{\tau},\frac{u_{\tau h}^{k+1}-u_{\tau h}^{k}}{\tau}\rangle + \frac12\langle\mathcal W_h'u_{\tau h}^{k+1},u_{\tau h}^{k}\rangle.
$$
This is an energy, {\em i.e.}, a positive quantity under the following \emph{Courant-Fridrichs-Lewy (CFL) condition} \cite{CoFrLe28PDMP} 
\begin{align}
  \big\langle \mathcal T' u_h,u_h\big\rangle 
  \geq \frac{\tau^2}{4}\big\langle{\cal W}_h' u_h,u_h \big\rangle
  \label{CFLold}
\end{align}
for any $u_h$ from the respective finite-dimensional subspace. The CFL typically bounds the time discretization step $\tau=\mathscr{O}(h_{\rm min})$ with $h_{\rm min}$ the smallest
element size on a FEM discretization. \mman{This method has frequently been
  used and analysed from various aspects, including comparison with implicit
  time discretizations, cf.\ e.g.\ \cite{KPCOP16TSDS,KPCG13GDAP}. However, }
the form of the discrete stored energy $\frac12\langle\mathcal W_h'u_{\tau h}^{k+1},u_{\tau h}^{k}\rangle$ makes this discretization
less suitable for the problem that we wish to consider in this paper where the stored energy is enhanced by some internal variables and
(possibly) nonlinear processes on them.


Therefore, we use another explicit discretization scheme, the so-called
\emph{leap-frog scheme}. To this end, we first write the velocity/stress
formulation of \eqref{IBVP1}. Introducing the velocity, $v=\DT u$  and 
  the stress tensor $\sigma:=\bbC e(u)$, we get 
\begin{subequations}\label{IBVP+}\begin{align}\label{IBVP1+}
    &&&  \varrho\DT v-{\rm div}\,\sigma=f\ \ \ \text{ and }\ \ \ 
    \DT\sigma=\bbC e(v)&&\text{on }\ \varOmega\ \text{ for }\
  t\in \mman{(}0,T],&&
  \\&&&\label{IBVP2+}
  \DT\sigma\vec{n}+\bbB 
  v
  =\DT g&&\text{on }\ \varGamma\ \text{ for }\
  t\in\mman{(}0,T],
  \\&&&\label{IBVP3+}
  v|_{t=0}=v_0,\ \ \ \sigma|_{t=0}=\sigma_0:=\bbC e(u_0)&&\text{on }\ \varOmega.
\end{align}\end{subequations}
In the abstract form \eqref{IVP}, when writing $\mathcal W=\mathscr{W}\circ
E$ with $E$ denoting the linear operator $u\mapsto(e,w):=(e(u),u|_\varGamma^{})$
and \mman{with $u|_\varGamma^{}$ denoting the trace of $u$ on the boundary
  $\varGamma$},
this reads as 
\begin{subequations}\label{IVP+}\begin{align}\label{IVP1+}
  &&&\mathcal T'\DT v+E^*\varSigma=F(t)&&\text{ for }\ t\in\mman{(}0,T],&&
  v|_{t=0}^{}=v_0,\ \ 
  \ \ \text{ and }\ \
  \\&&&\DT\varSigma=\mathscr{W}'Ev+\DT G(t)  
  &&\text{ for }\ t\in\mman{(}0,T],
  &&\varSigma|_{t=0}^{}=\varSigma_0:=\mathscr{W}'Eu_0,&&&&\label{IVP1++}
\end{align}\end{subequations}
where $E^*$ is the adjoint operator to $E$.
The stored energy
governing \eqref{IBVP+} is
$$
\mathscr{W}(e,w)=\int_\varOmega\frac12\bbC e{:}e\,\d x+
\int_\varGamma\frac12\bbB w{\cdot}w\,\d S
$$
while the external loading is now split into two parts acting differently,
namely
$$
\langle F(t),u\rangle=\int_\varOmega f(t)\cdot u\,\d x
\qquad\text{and}\qquad
\langle G(t),w\rangle=\int_\varGamma g(t)\cdot w\,\d x.
$$ 
Let us note that \eqref{IVP+} involves the
equation on $\varOmega$, as well as, the equation on $\varGamma$. Thus 
$\mathcal T$ is to be understood as the functional on $\varOmega\times\varGamma$,
that is trivial on $\varGamma$ since no \mman{inertia} is considered on
the $(d{-}1$)-dimensional boundary $\varGamma$. 
In particular, the ``generalized'' stress
$\varSigma=\mathscr{W}'Eu=(\bbC e(u),\bbB u|_\varGamma)$
contains, besides the bulk stress tensor, also the traction
stress vector. Relying on the linearity of $\mathscr{W}'$,
we have $\DT\varSigma=\mathscr{W}'Ev$ with $v=\DT u$, as used in
\eqref{IVP1++}.
\ccc
Let us note that the adjoint operator
$E^*:(\sigma,\varsigma)\mapsto \mathfrak{F}$ in
\eqref{IVP1+} with  the traction force $\varsigma=\bbB u|_\varGamma$
determines a bulk force $\mathfrak{F}$ as a functional on test
displacements $u$ by
\begin{align*}
  \int_\varOmega\mathfrak{F}\cdot u\,\d x
  =\left\langle \mathfrak{F},u\right\rangle
  &=\left\langle E^*(\sigma,\varsigma),u\right\rangle
  =\left\langle(\sigma,\varsigma),Eu\right\rangle
  =\left\langle(\sigma,\varsigma),
  (e(u),u|_\varGamma^{})\right\rangle
  \\[-.4em]&=\int_\varOmega\sigma:e(u)\,\d x+\int_\varGamma\varsigma\cdot u\,\d S
  =\int_\varGamma(\sigma\cdot\vec{n}+\varsigma)\cdot u\,\d S-\int_\varOmega{\rm div}\,\sigma\cdot u\,\d x\,,
\end{align*}
which clarifies the force $E^*\varSigma=\mathfrak{F}$ in \eqref{IVP1+}.
\eee
The leap-frog time discretization of
\eqref{IVP+} then reads as 
\begin{align}\label{elast-dyn-explicit+}
  \frac{\varSigma_{\tau h}^{k+1/2}-\varSigma_{\tau h}^{k-1/2}}\tau
  =\mathscr{W}_{\mman{h}}'E_hv_{\tau h}^{k}+D_{\tau h}^{k}
  \ \ \ \ \ \text{ and }\ \ \ \
  \TT'\frac{v_{\tau h}^{k+1}{-}v_{\tau h}^k}\tau+E_h^*\varSigma_{\tau h}^{k+1/2}
  =F_{\tau h}^{k+1/2},
\end{align}
where $\mathscr{W}_h$
and $E_h$ \cman{are} suitable FEM discretization\cman{s} \mman{of $\mathscr{W}$ and} $E$\mman{, cf.\
\eqref{suggestion-1} and \eqref{suggestion-2} below,}
and
\begin{align}\label{F-G}
  F_{\tau h}^{k+1/2}:=\frac1\tau\int_{k\tau}^{(k+1)\tau}\!\!\!\!\!\!\!F_h(t)\,\d t\
  \ \ \text{ and }\ \ \
  D_{\tau h}^{k}:=\frac1\tau\int_{(k-1/2)\tau}^{(k+1/2)\tau}\!\!\DT G_h(t)\,\d t
  =\frac{G_h((k{+}\frac12)\tau)-G_h((k{-}\frac12)\tau)}\tau\,.
\end{align}
\cman{As mentioned before, we} assume here that $v$ is discretized in an element-wise constant way so that
$\mathcal T$ leads to a diagonal matrix. In this case we do not need to employ numerical
integration to approximate the mass matrix. For higher-order discretizations, however,  
mass lumping is necessary so as to obtain explicit discretization schemes. 
We refer to \cite{BeJoTs02NFMF,JolTso08FEMD,Tsog99} for details in the case $G\equiv0$. 
The proposed FEM leads to a block diagonal matrix for $\mathcal W_h'=E_h^*\mathscr{W}_{\mman{h}}'E_h^{}$,
which means that the resulting scheme does not require the solution of a big linear system at each iteration in time. 
The spatial FEM discretization exploits regularity available in
linear elastodynamics, in particular that ${\rm div}\,\sigma$ and
$e(v)$ in \eqref{IBVP1+} live in $L^2$-spaces.
Moreover, the equations in \eqref{elast-dyn-explicit+} are decoupled in the
sense that, first, $\varSigma_{\tau h}^{k+1/2}$ is calculated from
the former equation and, second, $v_{\tau h}^{k+1}$ is calculated from the latter
equation assuming, that $(v_{\tau h}^k,\varSigma_{\tau h}^{k-1/2})$ is known
from the previous time step. For $k=0$, it starts from $v_{\tau h}^0=v_0$
and
from a half time step  
$\varSigma_{\tau h}^{1/2}=\varSigma_{\tau h}^{0}+\frac{\tau}{2} \mathscr{W}_{\mman{h}}'E_hv_{\tau h}^{0}$.
For the space discretization, the lower order $Q^{\rm div}_{k+1}-Q_{k}$ finite element is obtained for $k=0$ and in this case the velocity is discretized as piecewise
constant on rectangular or cubic
elements while the stress is discretized by piecewise bi-linear functions with
some continuities. Namely the normal component of the stress is continuous
across edges of adjacent elements while the tangential component is allowed
to be discontinuous. For more details about the space discretization we
refer the interested reader to \cite{BeJoTs02NFMF}. 
Alternative discretizations for the linear elasticity problem have been proposed by D.\,Arnold and his
collaborators who designed mixed finite elements for general
rectangular and triangular grids \cite{ArnoldWinther02, ArnoldAwanou05, ArnoldAwanouWinther08}. \mman{For tetrahedral leap-frog discretization of the elastodynamics see \cite{DelGli15AHOS}.}
In general, the leap-frog scheme has been frequently used in
geophysics to calculate seismic wave propagation with the finite differences
method, cf.\ e.g.\
\cite{Bohl02PVFD,Grav96SSWP,Viri84SHWP}. 

When taking the average (i.e.\ the sum with the weights $\frac12$ and
$\frac12$) of the second equation in \eqref{elast-dyn-explicit+}
in the level $k$ and $k{-}1$ tested by $v_{\tau h}^k$ and summing it with the
first equation in \eqref{elast-dyn-explicit+} tested by
$[\mathscr{W}_{\mman{h}}']^{-1}(\varSigma_{\tau h}^{k+1/2}{+}\varSigma_{\tau h}^{k-1/2})/2$,
we obtain
\begin{align*}&\frac1{2\tau}\big\langle[\mathscr{W}_{\mman{h}}']^{-1}\varSigma_{\tau h}^{k+1/2},
\varSigma_{\tau h}^{k+1/2}\big\rangle
-\frac1{2\tau}\big\langle[\mathscr{W}_{\mman{h}}']^{-1}\varSigma_{\tau h}^{k-1/2},
\varSigma_{\tau h}^{k-1/2}\big\rangle
\\&\qquad\qquad
=\Big\langle\frac{\varSigma_{\tau h}^{k+1/2}{+}\varSigma_{\tau h}^{k-1/2}}2,
E_h^{}v_{\tau h}^{k}\Big\rangle+\Big\langle
[\mathscr{W}_{\mman{h}}']^{-1}D_{\tau h}^k,
\frac{\varSigma_{\tau h}^{k+1/2}{+}\varSigma_{\tau h}^{k-1/2}}2\Big\rangle
\ \ \text{ and}
\\
&\Big\langle\TT'\frac{v_{\tau h}^{k+1}{-}v_{\tau h}^{k-1}}\tau,
v_{\tau h}^{k}\Big\rangle
+\Big\langle\frac{\varSigma_{\tau h}^{k+1/2}{+}\varSigma_{\tau h}^{k-1/2}}2,E_h^{}v_{\tau h}^{k}\Big\rangle
=\Big\langle\frac{ F_{\tau h}^{k+1/2}{+}F_{\tau h}^{k-1/2}}2,v_{\tau h}^{k}\Big\rangle\,.
\end{align*}
Summing it up, we get that the following discrete energy is conserved 
\begin{align}\label{leap-frog-energy}
  \frac12\langle
\ccc\TT\eee'
  v_{\tau h}^{k+1},v_{\tau h}^k\rangle+
\varPhi_h
(\varSigma_{\tau h}^{k+1/2})\ \ \text{ with }\ \ 
\varPhi_h(\varSigma)=\frac12\langle[\mathscr{W}_{\mman{h}}']^{-1}\varSigma,\varSigma\rangle\,.
\end{align}
Note that $\varPhi_h$ is the discrete stored energy expressed in terms of the
generalized stress. 
In contrast to  \eqref{elast-dyn-explicit}, this formulation allows for 
enhancement of the discrete stored energy by
some internal variables. 
The energy (\ref{leap-frog-energy}) is shown to be a positive quantity under the following CFL condition
\begin{align}
  \big\langle[ \mathcal W'_h ]^{-1}\varSigma_h,\varSigma_h\big\rangle
  \geq \frac{\tau^2}{4}\big\langle E^*_h \varSigma_h,(\TT')^{-1}  E^*_h \varSigma_h\big\rangle 
  \label{CFLold2}
\end{align}
for any $\varSigma_h$ from the respective finite-dimensional subspace.
Moreover, $F=0$ is often considered, which makes
the a-priori estimation easier.
Let us also note that the adjective ``leap-frog'' is sometimes
used also for the time-discretization \eqref{elast-dyn-explicit} if
written as a two-step scheme, cf.\ e.g.\ \cite[Sect.\,7.1.1.1]{CohPer17FEDG}.

The plan of this article is as follows: In Section~\ref{sec-int-var},
we complement the abstract system \eqref{IVP+} by another equation for
some internal variable 
and cast its weak formulation without relying on any regularity.
Then, in Section~\ref{sect-disc}, we extend the
two-step leap-frog discrete scheme \eqref{elast-dyn-explicit+}
to a suitable three-step scheme, and study the energy properties of the proposed scheme.
Then, in Section~\ref{sect-anal}, we   
prove the numerical stability of the 3-step staggered approximation scheme
and its convergence under the CFL condition 
\eqref{CFLnew}. 
Such an abstract scheme is then illustrated in Section~\ref{sec-exa}
on several examples from continuum mechanics, in particular
on models of plasticity, creep, diffusion,
and damage.
\mman{For illustration of computational efficiency, we refer to
  \cite{RoPaTs19ETDE} where this scheme was implemented
 for another problem, namely a delamination (i.e.\ interfacial damage).} 
It should be emphasized that, to the best of our knowledge, a rigorously
justified (as far as numerical stability and convergence) combination of the
explicit staggered discretization with nonlinear dissipative processes
on some internal variables is new, 
although occasionally some dissipative nonlinear phenomena can be found in
literature as in \cite{Scar04ETNP} for a unilateral contact, 
  in \cite{Bohl02PVFD} for a Maxwell viscoelastic rheology,
  in \cite{SePaPa18SEIF} for electroactive polymers,
  \mman{ in \cite{FaLeMa95MEIT} for an aeroelastic system,}
or in \cite{FePaFa01PACM} for general thermomechanical systems,
but without any numerical stability (a-priori estimates) and
convergence guaranteed. 



\section{Internal variables and their dissipative evolution.}\label{sec-int-var}

The concept of internal variables has a long tradition, \mman{cf.}\
\cite{Maug15SIVS}, and opens wide options for material modelling\mman{,
  cf.\ e.g.\ \cite{KruRou19MMCM,MieRou15RIST} and references therein}.
\mman{Typically, internal
  variables are governed by a parabolic-type}
1st-order evolution.
The \mman{abstract} system \eqref{IVP} is thus generalized to 
\begin{subequations}\label{IVP2}
  \begin{align}\label{IVP2Sigma}
  &\mathcal T'\DDT u+\mathcal W_u'(u,z)=F(t)&&
    \text{ for }\ t\in\mman{(}0,T],\ \ u|_{t=0}=u_0,\ \ \ \DT u|_{t=0}=v_0,
  \\&\label{IVP2z}\partial\varPsi(\DT z)+\mathcal W_z'(u,z)\ni0&&
  \text{ for }\ t\in\mman{(}0,T],\ \ z|_{t=0}=z_0.
\end{align}\end{subequations}
The inclusion in \eqref{IVP2z} refers to a possibility that the convex
(pseudo)potential $\varPsi$ of dissipative forces may be nonsmooth and then
its subdifferential $\partial\varPsi$ can be multivalued.

Combination of the 2nd-order evolution
\eqref{IVP} with such 1st-order evolution is to be \cman{handled} carefully.
In contrast to the implicit \mman{staggered} schemes, cf.\ \cite{RouPan17ECTD},
the constitutive equation is differentiated in time, cf.\ \eqref{IBVP1+},
and it seems necessary to use the staggered scheme so that
the internal-variable flow rule can be used without being
differentiated in time, even for a quadratic stored energy $\mathcal W$. 


Moreover, to imitate the leap-frog scheme, it seems suitable (or maybe
even necessary) that the stored energy $\mathcal W$ may be expressed
in terms of the generalized stress as
\begin{align}\label{ansatz}
  \mathcal W(u,z)=\varPhi(\varSigma,z)\ \ \text{ with }\ \ 
  \varSigma=\mathfrak C Eu,\ \text{ and }\ \ \varPhi(\cdot,z)\ \text{ and }\
  \varPhi(\varSigma,\cdot)\text{ quadratic},
\end{align}
where $\mathfrak C$ stands for a ``generalized'' elasticity
tensor and $E$ is an abstract gradient-type operator.
Typically $Eu=(e(u),u|_\varGamma)$ or simply $Eu=e(u)$ are considered here 
in the context of continuum mechanics at small strains,
cf.\ the examples in Sect.\,\ref{sec-exa}.
Here, $\varSigma$ may not directly enter the balance of forces and
is thus to be called rather as some ``proto-stress'',
while the actual generalized stress will be denoted by $S$.
For a relaxation of the last requirement of \eqref{ansatz} see Remark~\ref{rem-nonquadratic}
below.

Then, likewise \eqref{IVP+}, we can write the system \eqref{IVP2}
in the velocity/proto-stress formulation as
\begin{subequations}\label{IVP2+}
  \begin{align}\label{IVP2Sigma+}
&\DT\varSigma=\mathfrak CEv+\DT G(t)
    &&
    \text{ for }\ t\in\mman{(}0,T],\ \ &&
    \\ \label{IVP2S+}   &\mathcal T'\DT v+E^*S 
    =F(t)\ \ \text{ with }\ \
    S=\mathfrak C^*\varPhi_\varSigma'(\varSigma,z)\!\!\!\!\!\!&&
    \text{ for }\ t\in\mman{(}0,T],\ \ 
    \\&\label{IVP2z+}
    \partial\varPsi(\DT z)+\varPhi_z'(\varSigma,z)  
    \ni0&&
    \text{ for }\ t\in\mman{(}0,T],
    \\&\label{IVP2+IC}
    \varSigma|_{t=0}=\varSigma_0:=\mathfrak CEu_0+G(0)\,,\ \ \ v|_{t=0}=v_0\,,
    \ \ z|_{t=0}=z_0\,.\hspace*{-9em}
\end{align}\end{subequations}
Here $\varPhi_\varSigma'(\varSigma,z)$ is a ``generalized''
strain and, when multiplied by $\mathfrak C^*$, it becomes a 
generalized stress. \ccc Actually, \eqref{IVP+} is a special case of
\eqref{IVP2+} when $\varPhi=\varPhi(\varSigma)=
\int_\varOmega\frac12\bbC^{-1}\sigma:\sigma\,\d x+\int_\varGamma\bbB^{-1}\varsigma
\cdot\varsigma\,\d S$ and
$S=\varSigma=(\sigma,\varsigma)$ while $\mathfrak C=\mathscr{W}'=(\bbC,\bbB)$,
so that indeed $S=\mathfrak C^*\varPhi_\varSigma'(\varSigma)=
(\bbC,\bbB)(\bbC^{-1}\sigma,\bbB^{-1}\varsigma)=
(\sigma,\varsigma)=\varSigma$.
\eee

The energy properties of this system can be revealed by testing the particular
equations/inclusions in \eqref{IVP2+} by
$\varPhi_\varSigma'(\varSigma,z)$, $v$, and $\DT z$.
Thus, at least formally, we obtain
\begin{subequations}\label{energy-test}\begin{align}
&
\big\langle\varPhi_\varSigma'(\varSigma,z),\DT\varSigma
\big\rangle=\big\langle\varPhi_\varSigma'(\varSigma,z),\mathfrak CEv+\DT G\big\rangle
=\big\langle\mathfrak C^*\varPhi_\varSigma'(\varSigma,z),Ev\big\rangle
+\big\langle\varPhi_\varSigma'(\varSigma,z),\DT G\big\rangle,
\\  &
\ccc\big\langle\TT'\DT v,v\big\rangle\eee
+\big\langle\mathfrak C^*\varPhi_\varSigma'(\varSigma,z),Ev\big\rangle
  =\big\langle F(t),v\big\rangle\,,
  \\&\varXi(\DT z)
  +\big\langle\varPhi_z'(\varSigma,z),\DT z\big\rangle\le0\ \ \ \
  \text{ with }\ \ \ 
\varXi(\DT z):=\inf\big\langle\partial\varPsi(\DT z),\DT z\big\rangle\,.
\label{energy-test3}
\end{align}\end{subequations}
The functional $\varXi$ is a dissipative rate and 
the ``inf'' in it refers to the fact that the dissipative
potential $\varPsi$ can be nonsmooth and thus the 
subdifferential $\partial\varPsi$ can be multivalued even at $\DT z\ne0$,
otherwise an equality in \eqref{energy-test3} holds.
Summing it up and using
the calculus \ccc $\frac{\d}{\d t}\TT(v)=\langle\TT'v,\DT v\rangle
=\langle\TT'\DT v,v\rangle$ and \eee $\frac{\d}{\d t}\varPhi(\varSigma,z)=
\langle\varPhi_\varSigma'(\varSigma,z),\DT\varSigma\rangle
+\langle\varPhi_z'(\varSigma,z),\DT z\rangle$,
we obtain the following \ccc inequality \eee
for the energy, 
\begin{align}\label{energy}
  \frac{\d}{\d t}\big(\!\!\!\!\linesunder{\mathcal T(v)
    +\varPhi(\varSigma,z)}{kinetic and}{stored energies}\!\!\!\!\big)\
  +\!\!\!\!\!\!\linesunder{\varXi(\DT z)
  }{dissipation}{rate}\!\!\!\!\!\!
  \ccc\le\eee\!\!\linesunder{\big\langle F(t),v\big\rangle
   +\big\langle\varPhi_\varSigma'(\varSigma,z),\DT G\big\rangle}
  {power of}{external force}\,.
\end{align}
\ccc Actually, \eqref{energy-test3} and \eqref{energy} often hold as
equalities. \eee

Let us now formulate some abstract functional setting of the system
\eqref{IVP2+}. For some Banach spaces $\calS$, $\calZ$, and
$\calZ_1\supset\calZ$
and for a Hilbert space $\calH$, let 
$\varPhi:{\calS}\times\calZ\to\R$ be smooth and coercive, $\mathcal T:\calH\to\R$
be quadratic and coercive, and let $\Psi:\calZ\to[0,+\infty]$ 
be convex, lower semicontinuous, and coercive on $\calZ_1$, cf.\ \eqref{coerc+Lipschitz} below.
Intentionally, we do not want to rely on any regularity which is
usually at disposal in linear problems but might
be restrictive in some nonlinear problems. For this reason, we reconstruct
the abstract ``displacement'' and use \eqref{IVP2Sigma+} integrated in time,
i.e.
\begin{align}\label{IVP2Sigma++}
  \varSigma={\mathfrak C}Eu+G\ \ \ \text{ with }\ \
  u(t):=\int_0^t\!\!v(t)\,\d t+u_0\,.
\end{align}
Moreover, we still need another Banach space $\calE$
and define the Banach space $\calU:=\{u\in \calH;\
Eu\in \calE\}$ equipped with the standard graph norm. Then, by definition,
we have the continuous embedding $\calU\to\calH$ and the
continuous linear operator $E:\calU\to\calE$. We assume that
$\calU$ is embedded into $\calH$ densely, so that
$\calH^*\subset\calU^*$ and
that $\calH$ is identified with its dual $\calH^*$, so that
we have the so-called Gelfand triple
$$
\calU\subset\calH\cong\calH^*\subset\calU^*.
$$
We further consider the abstract elasticity tensor $\mathfrak C$
as a linear continuous operator $\calE\to\calS$.
Therefore $\mathfrak C Eu\in\calS$ provided $u\in \calU$ so that
the equation \eqref{IVP2Sigma++} is meant in $\calS$ and one needs
$G(t)\in\calS$.
Let us note that $\mathcal T':\calH\to\calH^*\cong\calH$,
$\varPhi_\varSigma':{\calS}\times\calZ\to{\calS}^*$,
$E^*:\calE^*\to\calU^*$, and
$\mathfrak C^*:\calS^*\to\calE^*$, so that
 $\mathcal T' v\in\calH^*$
provided $v\in\calH$
and also $S=\mathfrak C^*\varPhi_\varSigma'\in\calE^*$ and
$E^*S\in\calH^*$. In particular, the equation
\eqref{IVP2S+} can be meant in $\calH$ if integrated in time,
and one needs $F(t)$ valued in $\calH$.

\mman{For a Banach space $\calX$,} we will use the standard notation
$L^p(\mman{0,T;}\calX)$ for Bochner spaces of \mman{the} Bochner measurable
functions $\mman{[0,T]}\to\mathcal X$ whose norm is integrable with
the power $p$ or essentially bounded if $p=\infty$, and $W^{1,p}(\mman{0,T;}\calX)$
the space of functions from $L^p(\mman{0,T;}\calX)$ whose distributional time derivative
is also in  $L^p(\mman{0,T;}\calX)$. Also, $C^k(\mman{0,T;}\calX)$ will denote the space of
functions $\mman{[0,T]}\to\calX$ whose $k^{\rm th}$-derivative is continuous, and
$C_{\rm w}(\mman{0,T;}\calX)$ will denote the space of weakly continuous functions
$\mman{[0,T]}\to\mathcal X$. Later, we will also use ${\rm Lin}(\calU,\calE)$,
denoting the space of linear bounded operators $\calU\to\calE$
normed by the usual sup-norm.

A weak formulation of \eqref{IVP2S+} can be obtained after by-part integration
over the time interval $I=[0,T]$ when tested by a smooth function.
It is often useful to consider 
\begin{align}\label{z-ansatz}
  \varPhi(\varSigma,z)=\varPhi_0(\varSigma,z)+\varPhi_1(z)
  \ \ \ \text{ with }\ [\varPhi_0]_z':{\calS}\times\calZ
  \to\calZ_1^*\ \text{ and }\ \varPhi_1':\calZ\to\calZ^*
\end{align}
and to use integration by-parts for the term
$\langle\varPhi_1'(z),\DT z\rangle$. We thus arrive to the following definition. 

\begin{definition}[Weak solution to \eqref{IVP2+}.]\label{def}
  The quadruple $(u,\varSigma,v,z)\in C_{\rm w}(\mman{0,T;}{\calU})\times
  C_{\rm w}(\mman{0,T;}{\calS})\times C_{\rm w}(\mman{0,T;}{\calH})\times C_{\rm w}(\mman{0,T;}\calZ)$
with $\varPsi(\DT z)\in L^1(I)$ \mman{and $\DT z\in L^1(\mman{0,T;}\calZ_1)$} 
  will be called a weak solution to the initial-value problem
  \eqref{IVP2+} with \eqref{IVP2Sigma++}
  if $v=\DT u$ in the distributional sense, $\varSigma={\mathfrak C}Eu+G$
  holds a.e.\ on $I$, and if
 \begin{subequations}\label{weak}\begin{align}
&\label{weak-v}
\int_0^T\!\!\big\langle \varPhi_\varSigma'(\varSigma,z),\mathfrak C E\wt v\big\rangle_{{\calS}^*\times{\calS}}^{}
    -\big\langle\mathcal T'v,\DT{\wt v}\big\rangle_{{\calH}^*\times{\calH}}^{}
    \,\d t=\big\langle\mathcal T'v_0,\wt v(0)\big\rangle_{{\calH}^*\times{\calH}}^{}
    +\int_0^T\!\!\big\langle F,\wt v\big\rangle_{{\calH}^*\times{\calH}}^{}\,\d t
    \intertext{for any $\wt v\in C^1(\mman{0,T;}{\calH})\,\cap\, C(\mman{0,T;}{\calU})$ with $\wt v(T)=0$, and}
   &
   \int_0^T\!\!\varPsi(\wt z)+\big\langle[\varPhi_0]_z'(\varSigma,z),\wt z{-}
    \DT z\big\rangle_{\calZ_1^*\times\calZ_1}^{}\!
    +\big\langle\varPhi_1'(z),\wt z\big\rangle_{\calZ^*\times\calZ}^{}
    \,\d t+\varPhi_1(z_0)
    \ge\varPhi_1(z(T))+
    \!\int_0^T\!\!\varPsi(\DT z)\,\d t
    \label{weak-z}
 \end{align}\end{subequations}
 for any $\wt z\in C(\mman{0,T;}\calZ)$, where indices in the dualities
 $\langle\cdot,\cdot\rangle$ indicate the respective spaces in dualities, and
 if also $u(0)=u_0$, $\varSigma(0)=\varSigma_0$, and $z(0)=z_0$.
\end{definition}

Let us note that the remaining initial condition $v(0)=v_0$ is contained in
\eqref{weak-v}. Definition \ref{def} works successfully for $p>1$, i.e.\ for
rate-dependent evolution of the abstract internal variable $z$, so that
$\DT z\in L^p(\mman{0,T;}\calZ_1)$. For the rate-dependent evolution when $p=1$,
we would need to modify it. Here, we restrict ourselves to $p\ge2$,
because of the a-priori estimates in Proposition~\ref{prop1}.

\section{A three-step staggered time discretization}\label{sect-disc}

We derive in this section the leap-frog discretization of (\ref{IVP2+}a,b)
combined with an implicit discretization for \eqref{IVP2z+},
\mman{ using a fractional-step split (called also a staggered scheme)
  with a mid-point formula for \eqref{IVP2z+}}. 
Instead of a two-step scheme
\eqref{elast-dyn-explicit+}, we obtain a three-step scheme and
therefore, from now on, we abandon the convention
of a half-step notation as used in \eqref{elast-dyn-explicit+}
and write $k+1$ instead of $k+1/2$. 

\mman{To this aim, we consider sequences of nested
  finite-dimensional subspaces $S_h\subset\calS$,
  $V_h\subset\calH$, and $Z_h\subset\calZ$ where the values of the respective
  discrete variables \ccc $\varSigma_h$, $v_h$, and $z_h$ \eee will be, assuming that
  their unions are dense in the respective Banach spaces.
  We will use an
interpolation operator $I_h:{\rm Lin}(\calS,S_h)$ 
  and the
  embedding operator $J_h:Z_h\to\calZ$; it is important that the collection
  $\{J_h\}_{h>0}^{}$ is uniformly bounded and, since $\bigcup_{h>0}Z_h$ is dense in
  $\calZ$, the sequence $\{J_h\}_{h>0}^{}$ converges to the identity on $\calZ$
  strongly.} \mman{We consider $E_h\in{\rm Lin}(V_h,\calE)$.}
\cman{Let us note that we allow for  
a ``non-conformal'' approximation of $v$, i.e., $V_h\subset\calH$ is not necessarily a subspace of $\calU$. This is in agreement with discretizations of the velocity as in 
\cite{BeJoTs02NFMF,BeJoTs01MFES,BeRoTs09CRFD,ChLaQi15SDGM,Tsog99}.
}


Considering that we know from previous step
$\varSigma_{\tau h}^{\KKK},v_{\tau h}^k,z_{\tau h}^{k}$, then the proposed discretization scheme is
\begin{subequations}\label{suggestion}
\begin{align}\label{suggestion-1}
&\text{1) calculate $\varSigma_{\tau h}^{\KK}$:}&&\!\!\!
\frac{\varSigma_{\tau h}^{\KK}-\varSigma_{\tau h}^{\KKK}}\tau=
\mman{I_h}\mathfrak CE_h v_{\tau h}^{k}+D_{\tau h}^{k},
\\\label{suggestion-3}&\text{2) calculate $z_{\tau h}^{k+1}$:}&&\!\!\!
\mman{J_h^*}\partial\varPsi\Big(\frac{z_{\tau h}^{k+1}{-}z_{\tau h}^{k}}\tau\Big)+
\varPhi_z'\Big(\varSigma_{\tau h}^{\KK},\frac{z_{\tau h}^{k+1}{+}z_{\tau h}^{k}}{2}\Big)\ni0\,,
\\\nonumber
  &\text{3) calculate $v_{\tau h}^{k+1}$:}&&\!\!\! 
  \TT'\frac{v_{\tau h}^{k+1}{-}v_{\tau h}^{k}}{\tau}+E_h^\ast S_{\tau h}^{\KK}
  =F_{\tau h}^{\mman{k+1/2}}
  \ \text{ with }\ S_{\tau h}^{\KK}=\mathfrak C^*\mman{I_h^*}\varPhi_\varSigma'(\varSigma_{\tau h}^{\KK},\COLOR{z_{\tau h}^{k+1}}),&&
\label{suggestion-2}
\\&\text{\ \ \ \ \ \ \ \ \ and \ $u_{\tau h}^{k+1}$:}&&\!\!\!
u_{\tau h}^{k+1}=u_{\tau h}^{k}+\tau v_{\tau h}^{k+1},
\end{align}\end{subequations}
where $F_{\tau h}^{\mman{k+1/2}}$ and $D_{\tau h}^{k}$ are from \eqref{F-G}.
\cman{It seems important in the non-linear case to compute the variables in the order given above. We note however, that 
for  the linear viscoelastic problem with Maxwell rheology a scheme with a different ordering has been proposed 
  in \cite[Part~I, Sect.2]{Grob05MPOE}.} 
The potentials $\varPhi$, $\varPsi$, and $\TT$ are considered
  restricted on $S_h\times Z_h$, $Z_h$, and $V_h$, so that their corresponding 
  (sub)differentials $\varPhi_\varSigma'$, $\varPhi_z'$, $\partial\varPsi$, and
  $\TT'$ are valued in $S_h^*$, $Z_h^*$, and $V_h^*$, respectively.
  The particular equations/inclusion in \eqref{suggestion}
  are thus to be understood in $S_h$, $Z_h^*$, $V_h^*$,
\ccc $\calE^*$,  \eee
  and $V_h$,
  respectively.
  The only implicit equation is \eqref{suggestion-3}. \cman{Note \mman{however}
    that
  even this equation becomes explicit if there are no spatial gradients
  in $\varPhi$ and $\varPsi$}. In view of the definition of the
convex subdifferential, \eqref{suggestion-3} means the variational
inequality
\begin{align}\label{VI}
  \varPsi\big(\widetilde z\big)+\bigg\langle
\varPhi_z'\Big(\varSigma_{\tau h}^{\KK},\frac{z_{\tau h}^{k+1}{+}z_{\tau h}^{k}}{2}\Big),
  \widetilde z-\frac{z_{\tau h}^{k+1}{-}z_{\tau h}^{k}}\tau\bigg\rangle
  \ge\varPsi\Big(\frac{z_{\tau h}^{k+1}{-}z_{\tau h}^{k}}\tau\Big)
\end{align}
for any $\widetilde z\mman{\in V_h}$.
In any
case, the equation \eqref{suggestion-3} \mman{posseses a potential}
\begin{align}
z\mapsto\frac2\tau\varPhi\Big(\varSigma_{\tau h}^{\KK},\frac{z{+}z_{\tau h}^{k}}{2}\Big)+\varPsi\Big(\frac{z{-}z_{\tau h}^{k}}\tau\Big)
\end{align}
\mman{which is
  to be minimized on $Z_h$.}
Therefore the existence of a solution to \mman{this}
inclusion \eqref{suggestion-3}, or 
\mman{equivalently of the} variational inequality \eqref{VI},
can be shown by a direct method,
cf.\ also \cite{RouPan17ECTD}. 
\CHECK{The scheme \eqref{suggestion} is thus to be solved recurrently
  for $k=0,1,...,T/\tau-1$, starting from the initial conditions
  \eqref{IVP2+IC} assumed, for simplicity, to live in the \ccc respective
  finite-dimensional spaces.\eee}
  


The energy properties of this scheme can be obtained by imitating
\eqref{energy-test}--\eqref{energy}. More specifically, we proceed as 
follows: \mman{we} test \eqref{suggestion-1} by
$\frac12\varPhi_\varSigma'(\varSigma_{\tau h}^{\KK},\COLOR{z_{\tau h}^{k+1}})
+\frac12\varPhi_\varSigma'(\varSigma_{\tau h}^{\KKK},\COLOR{z_{\tau h}^{k}})$,
then \mman{test} \eqref{suggestion-3} by $\frac{z_{\tau h}^{k+1}-z_{\tau h}^{k}}\tau$, and
eventually \mman{test} the average of \eqref{suggestion-2} at the level $k{+}1$ and ${k}$
by $v_{\tau h}^{k}$. Using that $\varPhi(\cdot,z)$ and $\varPhi(\Sigma,\cdot)$
are quadratic as assumed in \eqref{ansatz}, we have
\begin{subequations}\label{test0}
  \begin{align}\nonumber
&\bigg\langle\frac{\varPhi_\varSigma'(\varSigma_{\tau h}^{\KK},\COLOR{z_{\tau h}^{k+1}})
  +\varPhi_\varSigma'(\varSigma_{\tau h}^{\KKK},\COLOR{z_{\tau h}^{k}})}2,
    \frac{\varSigma_{\tau h}^{\KK}\!-\varSigma_{\tau h}^{\KKK}}\tau\bigg\rangle
    \\&\nonumber\quad=
\bigg\langle\frac{\varPhi_\varSigma'(\varSigma_{\tau h}^{\KK},z_{\tau h}^{k})
  +\varPhi_\varSigma'(\varSigma_{\tau h}^{\KKK},z_{\tau h}^{k})}2,
    \frac{\varSigma_{\tau h}^{\KK}\!-\varSigma_{\tau h}^{\KKK}}\tau\bigg\rangle
    \\&\nonumber\qquad\qquad\qquad\qquad\qquad\qquad\qquad
    \COLOR{+}\frac\tau2\bigg\langle\frac{\varPhi_\varSigma'(\COLOR{\varSigma_{\tau h}^{\KK},z_{\tau h}^{k+1}})
  -\varPhi_\varSigma'(\COLOR{\varSigma_{\tau h}^{\KK},z_{\tau h}^{k}})}\tau,
    \frac{\varSigma_{\tau h}^{\KK}\!-\varSigma_{\tau h}^{\KKK}}\tau\bigg\rangle
    \\&
    \quad=
\frac{\varPhi(\varSigma_{\tau h}^{\KK},z_{\tau h}^{k})
  {-}\varPhi(\varSigma_{\tau h}^{\KKK},z_{\tau h}^{k})}\tau
  \COLOR{+}\frac\tau2\bigg\langle
  \frac{\varPhi_\varSigma'(\COLOR{\varSigma_{\tau h}^{\KK},z_{\tau h}^{k+1}})
  -\varPhi_\varSigma'(\COLOR{\varSigma_{\tau h}^{\KK},z_{\tau h}^{k}})}\tau,
    \frac{\varSigma_{\tau h}^{\KK}\!-\varSigma_{\tau h}^{\KKK}}\tau\bigg\rangle\,,
\label{test01}
\intertext{where we used also \eqref{suggestion-1}, and}
&\bigg\langle\varPhi_z'\Big(\varSigma_{\tau h}^{\KK},
\frac{z_{\tau h}^{k+1}{+}z_{\tau h}^{k}}{2}\Big),
\frac{z_{\tau h}^{k+1}-z_{\tau h}^{k}}\tau\bigg\rangle=
\frac{\varPhi(\varSigma_{\tau h}^{\KK},z_{\tau h}^{k+1})-
\varPhi(\varSigma_{\tau h}^{\KK},z_{\tau h}^{k})}\tau\,.
\end{align}\end{subequations}
Therefore, 
\mman{using again the particular equations/inclusion in 
\eqref{suggestion},
  we get respectively}
\begin{subequations}\label{test}
  \begin{align}\nonumber
&\frac{\varPhi(\varSigma_{\tau h}^{\KK},z_{\tau h}^{k})
      -\varPhi(\varSigma_{\tau h}^{\KKK},z_{\tau h}^{k})}\tau
    =
\bigg\langle\frac{\varPhi_\varSigma'(\varSigma_{\tau h}^{\KK},\COLOR{z_{\tau h}^{k+1}})
  +\varPhi_\varSigma'(\varSigma_{\tau h}^{\KKK},\COLOR{z_{\tau h}^{k}})}2,
\mman{I_h}\mathfrak CE_h v_{\tau h}^{k}+D_{\tau h}^k\bigg\rangle
\\[-.3em]&\hspace{14em}
\COLOR{-}\frac\tau2\bigg\langle\frac{\varPhi_\varSigma'(\COLOR{\varSigma_{\tau h}^{\KK},z_{\tau h}^{k+1}})
  -\varPhi_\varSigma'(\COLOR{\varSigma_{\tau h}^{\KK},z_{\tau h}^{k}})}\tau,
    \frac{\varSigma_{\tau h}^{\KK}\!-\varSigma_{\tau h}^{\KKK}}\tau\bigg\rangle
\label{test-1}
\,,
\\
&
\varXi
\Big(\frac{z_{\tau h}^{k+1}{-}z_{\tau h}^{k}}\tau\Big)
+\frac{\varPhi(\varSigma_{\tau h}^{\KK},z_{\tau h}^{k+1})-
\varPhi(\varSigma_{\tau h}^{\KK},z_{\tau h}^{k})}\tau\le0\,,\qquad\text{\mman{and}}
\label{test3}
\\
&\big\langle\TT'\frac{v_{\tau h}^{k+1}{-}v_{\tau h}^{k-1}}{2\tau},v_{\tau h}^{k}
\big\rangle+
\bigg\langle E_h^*\mathfrak C^*\mman{I_h^*}
\frac{\varPhi_\varSigma'(\varSigma_{\tau h}^{\KK},\COLOR{z_{\tau h}^{k+1}})
  {+}\varPhi_\varSigma'(\varSigma_{\tau h}^{\KKK},\COLOR{z_{\tau h}^{k}})}2,v_{\tau h}^{k}\bigg\rangle
=\big\langle F_{\tau h}^{\mman{k+1/2}},v_{\tau h}^{k}\big\rangle\,.
\label{test2}
\end{align}\end{subequations}
Let us also note that, if $\varPsi(0)=0$ is assumed, 
the substitution $\widetilde z=0$ into the inequality \eqref{VI}
gives $\varPsi(\frac{z_{\tau h}^{k+1}-z_{\tau h}^{k}}\tau)$
instead of the dissipation rate $\varXi(\frac{z_{\tau h}^{k+1}-z_{\tau h}^{k}}\tau)$
in \eqref{test3},
which is a suboptimal estimate except if $\varPsi$ is degree-1
positively homogeneous.

Summing \eqref{test} up, we benefit from the cancellation of the terms
$\pm\varPhi(\varSigma_{\tau h}^{\KK},z_{\tau h}^{k})$, which is the usual
\mman{goal and} attribute of 
\mman{well-designed} fractional-split schemes. Thus,
using also the simple algebra
\begin{align}\label{algebra-for-T}
  \big\langle\TT'(v_{\tau h}^{k+1}{-}v_{\tau h}^{k-1}),
v_{\tau h}^{k}\big\rangle=\big\langle\TT'v_{\tau h}^{k+1},v_{\tau h}^{k}\big\rangle
-\big\langle\TT'v_{\tau h}^{k},v_{\tau h}^{k-1}\big\rangle\,,
\end{align}
we obtain the analog of
\eqref{energy}, namely
\begin{align}\nonumber
&\frac{\langle\TT'v_{\tau h}^{k+1},v_{\tau h}^{k}\rangle
-\langle\TT'v_{\tau h}^{k},v_{\tau h}^{k-1}\rangle}{2\tau}
  +\frac{\varPhi(\varSigma_{\tau h}^{\KK},z_{\tau h}^{k+1})
    -\varPhi(\varSigma_{\tau h}^{\KKK},z_{\tau h}^{k})}\tau
+\varXi\Big(\frac{z_{\tau h}^{k+1}-z_{\tau h}^{k}}\tau\Big)
\\[-.3em]&\nonumber\qquad\qquad\qquad
\le\big\langle F_{\tau h}^{\mman{k+1/2}},v_{\tau h}^{k}\big\rangle
+
\bigg\langle\frac{\varPhi_\varSigma'(\varSigma_{\tau h}^{\KK},\COLOR{z_{\tau h}^{k+1}})
  +\varPhi_\varSigma'(\varSigma_{\tau h}^{\KKK},z_{\tau h}^{k})}2,D_{\tau h}^k\bigg\rangle
\\[-.1em]&\qquad\qquad\qquad\quad\
\COLOR{-}\frac\tau2\bigg\langle\frac{\varPhi_\varSigma'(\COLOR{\varSigma_{\tau h}^{\KK},z_{\tau h}^{k+1}})
  -\varPhi_\varSigma'(\COLOR{\varSigma_{\tau h}^{\KK},z_{\tau h}^{k}})}\tau,
    \frac{\varSigma_{\tau h}^{\KK}\!-\varSigma_{\tau h}^{\KKK}}\tau\bigg\rangle\,.
 \label{energy-disc}\end{align}
If $\varPsi$ is smooth except possibly at zero, there is even equality
in \eqref{energy-disc}.

Considering some approximate values $\{z_{\tau h}^k\}_{k=0,...,K}$ of the variable
$z$
with $K=T/\tau$, we define the piecewise-constant and the piecewise affine 
interpolants respectively by
\begin{subequations}\label{def-of-interpolants}
\begin{align}\label{def-of-interpolants-}
&&&
\overline{z}_{\tau h}(t)= z_{\tau h}^k,\qquad\ \
\underline{\overline z}_{\tau h}(t)=\frac12z_{\tau h}^k+\frac12z_{\tau h}^{k-1},
\qquad\ \text{ and}
&&
\\&&&\label{def-of-interpolants+}
z_{\tau h}(t)=\frac{t-(k{-}1)\tau}\tau z_{\tau h}^k
+\frac{k\tau-t}\tau z_{\tau h}^{k-1}
&&\hspace*{-7em}\text{for }(k{-}1)\tau<t\le k\tau.
\end{align}\end{subequations}
Similar meaning is implied for $\varSigma_{\tau h}$, $v_{\tau h}$,
$\overline\varSigma_{\tau h}$, $\overline v_{\tau h}$, $\overline F_{\tau h}$, etc. 
%
The discrete scheme \eqref{suggestion} can be written in a ``compact'' form as 
\begin{subequations}\label{suggestion+}
\begin{align}\label{suggestion-1+}
  &
  \DT\varSigma_{\tau h}=\mman{I_h}
  \mathfrak CE_h \overline v_{\tau h}+\DT G_{\tau h}\ \ \text{ and }\ \
  \DT u_{\tau h}=\overline v_{\tau h},
\\\label{suggestion-3+}
&
\mman{J_h^*}\partial\varPsi\big(\DT z_{\tau h}\big)+
\varPhi_z'\big(\overline\varSigma_{\tau h},\underline{\overline z}_{\tau h}\big)\ni0\,,
\\\label{suggestion-2+}
  &
  \TT'\DT v_{\tau h}+E_h^\ast\overline S_{\tau h}
  = \cman{{\overline F}_{\tau h}} 
  \ \ \text{ with }\ \ \overline S_{\tau h}
  =\mathfrak C^*\mman{I_h^*}\varPhi_\varSigma'(\overline\varSigma_{\tau h},\COLOR{\overline z_{\tau h}})\,
\end{align}\end{subequations}
\mman{to be valid a.e.\ on the time interval $[0,T]$.}

\section{Numerical stability and convergence}
\label{sect-anal}

Because the 
energy \eqref{leap-frog-energy} involves now also the internal variable, 
the CFL condition 
\mman{has to be modified. More specifically, we assume that}
\begin{align}
\exists\,\eta>0\ \forall\,\varSigma_h\mman{\in S_h},\ z_h\mman{\in Z_h}:\ \ 
  \varPhi(\varSigma_h,z_h)\ge\frac{\tau^2}{4{-}\eta}
      \big\langle E^*_h S_h,(\TT')^{-1}  E^*_h S_h\big\rangle_{{\calH}^*\times{\calH}}^{}
      \ \ \text{ with }\ S_h=\mathfrak C^*\mman{I_h^*}\varPhi'_{\Sigma} (\varSigma_h,z_h),
  \label{CFLnew}
\end{align}
where
$\varSigma_h$ and $z_h$ are considered
from the corresponding finite-dimensional subspaces. 
Let us still introduce the Banach space
${\mathcal X}:=\{X\in{\calS}^*;\ E^*{\mathfrak C}^*X\in{\calH}^*\}$.
We further assume
$\mathfrak C\in{\rm Lin}({\calE},{\calS})$
invertible \mman{and that
the collection of the interpolation operators $\{I_h:\calS\to S_h\}_{h>0}$ is
bounded in ${\rm Lin}(\calS,\calS)$.}

\begin{proposition}[Numerical stability.]\label{prop1}
  Let $F$ be constant in time, valued in ${\calH}^*$,
  $G\!\in\!W^{1,1}(\mman{0,T;}{\calS})$,
  $u_0\in\mathcal\calU$ so that $\varSigma_0={\mathfrak C}Eu_0\in
  {\calS}$, $v_0\in\calH$, $z_0\in\calZ$, the functionals
  $\mathcal T$, $\varPhi$, and $\varPsi$ be coercive and 
  $\varPhi_\varSigma'(\varSigma,\cdot)$ be Lipschitz continuous uniformly for
  $\varSigma\in\calS$ in the sense
\begin{subequations}\label{coerc+Lipschitz}\begin{align}\nonumber
    &&&\exists\epsilon>0\ p\ge2\ \forall (\varSigma,v,z)\in
    {\calS}{\times}{\calH}{\times}\calZ:
    \\&&&&&\hspace*{-9em}\mathcal T(v)\ge\epsilon\|v\|_{\calH}^2,\ \
    \varPhi(\varSigma,z)\ge\epsilon\|\varSigma\|_{\calS}^2+\epsilon\|z\|_\calZ^2,
    \ \
  \varPsi(z)\ge\epsilon\|z\|_{\calZ_1}^p,&&&&
  \\&&&\label{ass-Phi_S}
  \exists C\ \forall \varSigma\in\calS,\ z\in \calZ:&&
  \big\|\varPhi_\varSigma'(\varSigma,z)\big\|_{{\calS}^*}^{}\le
  C\big(1+\|\varSigma\|_{\calS}^{}+\|z\|_\calZ^{}\big),
  \\&&&\exists\ell\in\R\ \forall \varSigma\in\calS,\ z,\wt z\in \calZ:&&
    \big\|\varPhi_\varSigma'(\varSigma,z)-\varPhi_\varSigma'(\varSigma,\wt z)\big\|_{{\calS}^*}^{}
    \le\ell\|z-\wt z\big\|_{\calZ_1}^{}.
    \label{Lipschitz}
  \end{align}\end{subequations}
  Let also the CFL condition \eqref{CFLnew} hold with
  $\tau>0$ sufficiently small (in order to make the discrete Gronwall
  inequality effective).
Then the following a-priori estimates hold:
\begin{subequations}\label{est}
  \begin{align}\label{est0}
    &\|u_{\tau h}\|_{W^{1,\infty}(\mman{0,T;}\mathcal H)}^{}\le C\,,
    \\\label{est1}
    &\|\varSigma_{\tau h}\|_{L^\infty(\mman{0,T;}\calS)}^{}\le C\ \ \text{ and }\ \
    \|\DT\varSigma_{\tau h}\|_{L^1(\mman{0,T;}\mathcal X^*)}^{}\le C,
    \\&\label{est2}\|v_{\tau h}\|_{L^\infty(\mman{0,T;}\calH)}^{}\le C\ \ \text{ and }\ \
    \|\mathcal T'\DT v_{\tau h}\|_{L^\infty(\mman{0,T;}\calU^*)}^{}\le C,
    \\&\label{est3}\|z_{\tau h}\|_{L^\infty(\mman{0,T;}\calZ)}^{}\le C
    \ \ \text{ and }\ \
    \|\DT z_{\tau h}\|_{L^p(\mman{0,T;}\calZ_1)}^{}\le C.
\end{align}\end{subequations}    
\end{proposition}

\noindent{\it Proof.}
The energy 
{imbalance} that we have here is \eqref{energy-disc} which can be re-written as
 \begin{align}\nonumber
&\frac{\mathfrak E_h^{\KK}{-}\mathfrak E_h^{\KKK}}{\tau}+\varXi
\Big(\frac{z_{\tau h}^{k+1}{-}z_{\tau h}^{k}}\tau\Big)
\le\big\langle F_{\tau h}^{k},v_{\tau h}^{k}\big\rangle_{{\calH}^*\times{\calH}}^{}
+\bigg\langle\frac{\varPhi_\varSigma'(\varSigma_{\tau h}^{\KK},\COLOR{z_{\tau h}^{k+1}})
  +\varPhi_\varSigma'(\varSigma_{\tau h}^{\KKK},z_{\tau h}^{k})}2,D_{\tau h}^k\bigg\rangle_{{\calS}^*\times{\calS}}^{}
\\[-.1em]&\qquad\qquad\qquad\qquad\qquad\qquad
\COLOR{-}\frac\tau2\bigg\langle\frac{\varPhi_\varSigma'(\COLOR{\varSigma_{\tau h}^{\KK},z_{\tau h}^{k+1}})
  -\varPhi_\varSigma'(\COLOR{\varSigma_{\tau h}^{\KK},z_{\tau h}^{k}})}\tau,
    \frac{\varSigma_{\tau h}^{\KK}\!-\varSigma_{\tau h}^{\KKK}}\tau
   \bigg\rangle_{{\calS}^*\times{\calS}}^{}
 \label{en1}\end{align}
 with an analog of the energy \eqref{leap-frog-energy}, namely
 \begin{align}
   \mathfrak E_h^{\KK}  =
   \frac12\langle\TT'v_{\tau h}^{k+1},v_{\tau h}^{k}\rangle_{{\calH}^*\times{\calH}}^{}
 +\varPhi(\varSigma_{\tau h}^{\KK},z_{\tau h}^{k+1}).
  \label{en2}
\end{align}
 We need to show that $\mathfrak E_h^{k+1} $ is indeed 
 {a sum of the kinetic and the stored energies}
 { at least up to some positive coefficients}. To do so, like e.g.\
\cite[Lemma~4.2]{Scar04ETNP} or \cite[Sect.~6.1.6]{Tsog99}, let us write
\begin{align}
  \langle\TT'v_{\tau h}^{k+1},v_{\tau h}^{k}\rangle 
  &=  \Big\langle\TT' \frac{v_{\tau h}^{k+1}+v_{\tau h}^{k} }{2} ,\frac{v_{\tau h}^{k+1}+v_{\tau h}^{k} }{2}\Big\rangle
  -  \Big\langle\TT' \frac{v_{\tau h}^{k+1}-v_{\tau h}^{k} }{2} ,\frac{v_{\tau h}^{k+1}-v_{\tau h}^{k} }{2}\Big\rangle
  \nonumber \\\nonumber
  &=  \Big\langle \TT' \frac{v_{\tau h}^{k+1}+v_{\tau h}^{k} }{2} ,\frac{v_{\tau h}^{k+1}+v_{\tau h}^{k} }{2}\Big\rangle
  \\
  &\ - \frac{\tau^2}{4}\big\langle  E^*_h\big(\mathfrak C^*\mman{I_h^*}\varPhi'_\Sigma(\varSigma^{\KK}_{\tau h},\COLOR{z_{\tau h}^{k+1}}){-}F_{\tau h}^{\mman{k+1/2}}\big),
  (\TT')^{-1} E^*_h\big(\mathfrak C^*\mman{I_h^*}\varPhi'_\Sigma( \varSigma^{\KK}_{\tau h},\COLOR{z_{\tau h}^{k+1}}){-}F_{\tau h}^{\mman{k+1/2}}\big)\big\rangle\,,
\label{en3}
\end{align}
where all the duality pairings are between ${\calH}^*$ and
${\calH}$; 
here also \eqref{suggestion-2} has been used.
Thus, using also $\TT(v)=\frac12\langle\TT'v,v\rangle$,
we can write the energy \eqref{en2} as
\begin{align}\nonumber
  \mathfrak E_h^{\KK}&=
\TT(v_{\tau h}^{k+1/2})
+a_{\tau h}^{\KK}\varPhi(\varSigma_{\tau h}^{\KK},z_{\tau h}^{k+1})
+\frac{\tau^2}{2}\big\langle(\TT')^{-1}E^*_h\mathfrak C^*\mman{I_h^*}
\varPhi'_\Sigma(\varSigma^{\KK}_{\tau h},\COLOR{z_{\tau h}^{k+1}}),F_{\tau h}^{\mman{k+1/2}}\big\rangle
- \frac{\tau^2}{4}\|F_{\tau h}^{\mman{k+1/2}}\|_{\calH}^2   \\
  & \text{ with }\ \ a_{\tau h}^{\KK}:=
  1-\frac{\tau^2}{4}  \frac{  \big\langle  E^*_h \mathfrak C^*\mman{I_h^*}\varPhi'_\Sigma( \varSigma^{\KK}_{\tau h},\COLOR{z_{\tau h}^{k+1}}), (\TT')^{-1} E^*_h \mathfrak C^* \mman{I_h^*}\varPhi'_\Sigma( \varSigma^{\KK}_{\tau h},\COLOR{z_{\tau h}^{k+1}})\big\rangle}{\varPhi(\varSigma_{\tau h}^{\KK},z_{\tau h}^{k+1}) }\ge\eta
  \label{en4}\end{align}
and with $v_{\tau h}^{k+1/2}:=\frac12v_{\tau h}^{k+1}+\frac12v_{\tau h}^{k}$.
The energy $\mathfrak E_h^{\KK}$
        {yields a-priori estimates} if the
coefficient $a_{\tau h}^k$ 
is non-negative, which is just ensured by our CFL
  condition
%
  \eqref{CFLnew} used for $\varSigma_h=\varSigma_{\tau h}^{\KK}$ \mman{and}
  $z_h=z_{\tau h}^{k+1}$. 
  Here $\eta>0$ \mman{in \eqref{en4}} is just from \eqref{CFLnew}. 


  Altogether, summing \eqref{en1} for $k=0,...,l
  \le T/\tau-1$ and using \eqref{en4}, we obtain the estimate
 \begin{align}\nonumber
   &
   \epsilon\bigg(
   \big\|v_{\tau h}^{\mman{l+1/2}}\big\|_{\calH}^2
   +a_{\tau h}^{\mman{l+1}}\big\|\varSigma_{\tau h}^{\mman{l+1}}\big\|_{\calS}^2
   +a_{\tau h}^{\mman{l+1}}\big\|z_{\tau h}^{\mman{l+1}}\big\|_\calZ^2
   +\tau\sum_{k=0}^{\mman{l}}
   \Big\|\frac{z_{\tau h}^{k+1}{-}z_{\tau h}^{k}}\tau\Big\|_{\calZ_1}^p\bigg)
\\&\nonumber\quad
\le
\frac{\tau^2}{4}\|F_{\tau h}^{\mman{l+1/2}}\|_{\calH}^2
-\frac{\tau^2}{2}\big\langle(\TT')^{-1}E^*_h\mathfrak C^*\mman{I_h^*}
\varPhi'_\Sigma(\varSigma_{\tau h}^{\mman{l+1}},\COLOR{z_{\tau h}^{\mman{l+1}}}),F_{\tau h}^{\mman{l+1/2}}\big\rangle
-\frac{\tau^2}{2}\big\langle(\TT')^{-1}E^*_h\mathfrak C^*\mman{I_h^*}
\varPhi'_\Sigma(\varSigma^{0}_{\tau h},\COLOR{z_{\tau h}^{0}}),F_{\tau h}^{\mman{1/2}}\big\rangle
\\[-.5em]&\nonumber\quad\
   +\TT(v_{\tau h}^{-1/2})
   +a_{\tau h}^{0}\varPhi(\varSigma_{\tau h}^{0},z_{\tau h}^{0})
+\tau\sum_{k=0}^{\mman{l}}\bigg(\big\langle F_{\tau h}^{\mman{k+1/2}},v_{\tau h}^{k}\big\rangle
\\[-.5em]&
\quad\
   +\frac12\big\|\varPhi_\varSigma'(\varSigma_{\tau h}^{\KK},\COLOR{z_{\tau h}^{k+1}})
  +\varPhi_\varSigma'(\varSigma_{\tau h}^{\KKK},z_{\tau h}^{k})\big\|_{{\calS}^*}
\big\|D_{\tau h}^k\big\|_{{\calS}}^{}
+\frac\tau2\ell\Big\|\frac{z_{\tau h}^\mman{k+1}{-}z_{\tau h}^\mman{k}}\tau\Big\|_{\calZ_1}
\Big\|
\frac{\varSigma_{\tau h}^{\KK}\!{-}\varSigma_{\tau h}^{\KKK}}\tau\Big\|_{{\calS}}^{}\bigg)\,,
 \label{overall}\end{align}
 where $\epsilon$, $p$, $\ell$ and $a_{\tau h}^{\mman{l+1}}$ come from
 \eqref{coerc+Lipschitz} and \eqref{en4}. Here we also have used that
   the collection $\{I_h\}_{h>0}$ is bounded.
Using \eqref{ass-Phi_S}, we estimate 
$\|\varPhi_\varSigma'(\varSigma_{\tau h}^{\KKK},z_{\tau h}^{k})\|_{{\calS}^*}
\le C(1+\|\varSigma_{\tau h}^{\KKK}\|_{\calS}^2
+\|z_{\tau h}^{\KKK}\|_\calZ^2)$
\mman{and $\|\varPhi_\varSigma'(\varSigma_{\tau h}^{\KK},z_{\tau h}^{k\cman{+1}})
\|_{{\calS}^*}
\le C(1+\|\varSigma_{\tau h}^{\KK}\|_{\calS}^2
+\|z_{\tau h}^{\KK}\|_\calZ^2)$,} and then use the summability of
$\|D_{\tau h}^k\big\|_{{\calS}}^{}$ needed for the discrete Gronwall inequality;
here the assumption $\DT G\in L^1(\mman{0,T;}{\calS})$ is needed. The last
term in \eqref{overall} is to be estimated by the \REPLACE{H\"older}{Young}
inequality as 
\begin{align*}
\frac\tau2\ell\Big\|\frac{z_{\tau h}^\mman{k+1}{-}z_{\tau h}^\mman{k}}\tau\Big\|_{\calZ_1}
\Big\|
\frac{\varSigma_{\tau h}^{\KK}\!{-}\varSigma_{\tau h}^{\KKK}}\tau\Big\|_{\calS}^{}
&\mman{\le\frac\epsilon2\Big\|\frac{z_{\tau h}^\mman{k+1}{-}z_{\tau h}^\mman{k}}\tau\Big\|_{\calZ_1}^2
+\frac{\ell^2}{8\epsilon}\big\|\varSigma_{\tau h}^{\KK}\!{-}\varSigma_{\tau h}^{\KKK}\big\|_{\calS}^2}
\\&\le\frac\epsilon2\Big\|\frac{z_{\tau h}^\mman{k+1}{-}z_{\tau h}^\mman{k}}\tau\Big\|_{\calZ_1}^p
+C_{p,\epsilon,\ell}^{}\big(1+
\|\varSigma_{\tau h}^{\KK}\|_{\calS}^2
+\|\varSigma_{\tau h}^{\KKK}\|_{\calS}^2\big)
\end{align*}
with some $C_{p,\epsilon,\ell}^{}$ depending on $p$, $\epsilon$, and $\ell$. Here
we needed $p\ge2$; note that this is related with the specific explicit time
discretization due to the last term in \eqref{energy-disc} but not with the
problem itself. Then we use
the discrete Gronwall inequality
to obtain the former estimates in (\ref{est}b,c) and the estimates
(\ref{est}a,d).
Using the discrete Gronwall inequality is a bit
tricky because of the term $\|v_{\tau h}^{l\mman{+}1/2}\|_{\calH}^2$ \COMMENT{Now, $\|v_{\tau h}^{l+1/2}\|_{\calH}^2$ appears on the left-hand side of \eqref{overall} -- OK, CHANGED - THANKS}
on the left-hand side of \eqref{overall} while there is
$v_{\tau h}^{k}$ instead of $v_{\tau h}^{k{+}1/2}$ \COMMENT{now, it should read "instead of $\|v_{\tau h}^{k+1/2}\big\|_{\calH}^2$" -- SO I CHANGED $v_{\tau h}^{k{-}1/2}$ TO $v_{\tau h}^{k{+}1/2}$ BUT DID NOT WRITE $\|v_{\tau h}^{k+1/2}\big\|_{\calH}^2$ BECAUSE SUCH TERM IS NOT IN \eqref{overall}, OK?}
\cman{on the right-hand side of
  \eqref{overall}}. To cope with it,
we \mman{have to} rely on $F$ being constant (as assumed). We prove the
estimate for $l=1$, then we sum up \eqref{overall} for $l{+}1$ and $l$ to get 
$\big\langle F_{\tau h}^{\mman{k+1/2}},v_{\tau h}^{k\mman{+}1/2}\big\rangle$ \COMMENT{now, it should be $\big\langle F_{\tau h}^{\mman{k+1/2}},v_{\tau h}^{k+1/2}\big\rangle$  -- OK, CHANGED - THANKS} also
on the right-hand side. 
Note also that, in view of \eqref{algebra-for-T} for $k=0$, we have
  obtained the term $\TT(v_{\tau h}^{-1/2})$ on the right-hand side of
  \eqref{overall} which, however, can simply be ignored if taking the
  ``fictitious''
  velocity at level $k=-1$ as $-v_{\tau h}^{0}=-v_h^0$.

The equation $\DT\varSigma_{\tau h}=\mman{I_h}\mathfrak CE_h\overline v_{\tau h}
+\DT G_{\tau h}$ gives the latter estimate in
\eqref{est1} by estimating
\begin{align}\nonumber
  \int_0^T\!\!
  \big\langle\DT\varSigma_{\tau h},X\big\rangle_{\mathcal X^*\times\mathcal X}\,\d t
&=\int_0^T\!\!\big\langle\mman{I_h}\mathfrak CE_h\overline v_{\tau h}
+\DT G_{\tau h},X\big\rangle_{\mathcal X^*\times\mathcal X}\,\d t
\\&=\int_0^T\!\!\big\langle\overline v_{\tau h},E_h^*\mathfrak C^*\mman{I_h^*}X\big\rangle_{\calH\times\calH^*}\,\d t
+\int_0^T\!\!\big\langle\DT G_{\tau h},X\big\rangle_{\mathcal X^*\times\mathcal X}\,\d t
\label{dual-est1}\end{align}
for $X\in L^\infty(\mman{0,T;}\mathcal X)$ and using also the already proved boundedness
of $\overline v_{\tau h}$ in $L^\infty(\mman{0,T;}\calH)$ and the assumed
boundedness of $E_h$ uniform in $h>0$;
here we used also that $\DT\varSigma_{\tau h}(t)\in\calS\subset\mathcal X^*$.

Eventually, the already obtained estimates \eqref{ass-Phi_S} 
give $\varPhi_\varSigma'(\overline\varSigma_{\tau h},\COLOR{\overline z_{\tau h}})$
bounded in $L^\infty(\mman{0,T;}\calS^*)$. Therefore 
$\overline S_{\tau h}=\mathfrak C^*\mman{I_h^*}\varPhi_\varSigma'(\overline\varSigma_{\tau h},\COLOR{\overline z_{\tau h}})$ is bounded in $L^\infty(\mman{0,T;}\calE^*)$, hence
$E_h^*\overline S_{\tau h}$ is bounded in $L^\infty(\mman{0,T;}\calU^*)$,
so that $\mathcal T'\DT v_{\tau h}=\overline F_{\tau h}-E_h^*\overline S_{\tau h}$
gives the latter estimate in \eqref{est2}.
$\hfill\Box$

\bigskip

We will \mman{now} need the approximation propert\mman{ies for $h\to0$:} 
\begin{subequations}\begin{align}\label{Eh->E}
  &v\in\calU,\ \,v_h\in V_h,\ \ \,v_h\to v\ \ \text{ in }\ \calH\ \ \
  \Rightarrow\ \ \ 
  E_hv_h\to Ev\ \ \text{ in }\ \calE\ \ \text{ with }\
  E\in{\rm Lin}(\calU,\calE)\,,
  \\&\mman{\Sigma\in\calS,\ \Sigma_h\in S_h,\ \Sigma_h\to\Sigma\ \
    \text{ in }\ \calS\ \ \ \Rightarrow\ \ \ 
I_h\Sigma_h\to\Sigma\ \ \ \;\text{ in }\ \calS\,.}
\end{align}\end{subequations}

\begin{proposition}[Convergence.]\label{prop2}
  Let \eqref{z-ansatz} and \eqref{Eh->E} hold, all the involved Banach spaces
  be separable, 
  and the assumptions of Proposition~\ref{prop1} hold. Moreover, let
\begin{subequations}\label{ass-Phi'}\begin{align}\nonumber
    &\forall z\In\calZ:\ \varPhi_\varSigma'(\cdot,z)\ \text{continuous linear, and }
    \varPhi_\varSigma'
    :\calS{\times}\calZ_0\to{\rm Lin}(\calS,\calS^*)\ \text{continuous}
  \\&\label{ass-Phi-Sigma}\qquad\qquad\qquad\qquad\qquad
\text{or }\ \varPhi_\varSigma':\calS\times\calZ\to\calS^*\
  \text{ is continuous linear,}
  \\\nonumber
  &\forall z\In\calZ:\ [\varPhi_0]_z'(\cdot,z)\ \text{continuous linear, and }
      [\varPhi_0]_z'
      :\calS{\times}\calZ_0\to{\rm Lin}(\calZ_0,\calZ_1)\ \text{continuous}
  \\&\qquad\qquad\qquad\qquad\qquad\text{or }\ [\varPhi_0]_z':\calS\times\calZ\to\calZ_1^*\
  \text{  is continuous linear, and}
  \label{ass-Phi0-z}
  \\&\varPhi_1':\calZ\to\calZ^*\ \text{ is linear continuous},
  \end{align}\end{subequations}
for some Banach space $\calZ_0$ into which $\calZ$ is embedded compactly,
where $\varPhi_0$ and $\varPhi_1$ are from \eqref{z-ansatz}.
Then there is a selected subsequence, again denoted
  $\{(u_{\tau h},\varSigma_{\tau h},v_{\tau h},z_{\tau h})\}_{\tau>0}$
  converging weakly* in the topologies indicated in the estimates \eqref{est}
  to some $(u,\varSigma,v,z)$.
  Moreover, any $(u,\varSigma,v,z)$ obtained as such a limit is a weak
  solution according Definition~\ref{def}.
  \end{proposition}

\noindent{\it Proof.}
By the Banach selection principle, we can select the weakly* converging
subsequence as claimed; here the separability of the involved Banach spaces 
is used.

Referring to the compact embedding $\calZ\subset\calZ_0$ used in
the former option in (\ref{ass-Phi'}a,b) and relying on a generalization
the Aubin-Lions compact-embedding theorem with
$\COLOR{\DT{\overline z}_{\tau h}}$
being bounded in the space of the $\calZ_1$-valued measures on $I$,
cf.\ \cite[Corollary 7.9]{Roub13NPDE}, we have $\COLOR{\overline z_{\tau h}}\to z$
strongly in $L^{r}(\mman{0,T;}\calZ_1)$ for any $1\le r<+\infty$. 

Further, we realize that the approximate solution satisfy identities/inequality
analogous to what is used in Definition~\ref{def}.
In view of \eqref{weak-v}, the equations \eqref{suggestion-2+} now
  means
  \begin{subequations}\label{weak+}\begin{align}
&\label{weak-v-disc}
\int_0^T\!\!\big\langle\varPhi_\varSigma'(\overline\varSigma_{\tau h},\COLOR{\overline z_{\tau h}}),\mman{I_h}\mathfrak C E_h\wt v\big\rangle_{{\calS}^*\times{\calS}}^{}
    -\big\langle\mathcal T'v_{\tau h},\DT{\wt v}\big\rangle_{{\calH}^*\times{\calH}}^{}
    \,\d t=\big\langle\mathcal T'v_0,\wt v(0)\big\rangle_{{\calH}^*\times{\calH}}^{}
    +\int_0^T\!\!\big\langle F_{h},\wt v\big\rangle_{{\calH}^*\times{\calH}}^{}\,\d t
    \intertext{for any $\wt v\in C^1(\mman{0,T;}{\calH})
      $ valued in $V_h$ 
      and with $\wt v(T)=0$.
\REPLACE{In view of}{Like in} \eqref{weak-z}, the inclusion
      \eqref{suggestion-3+} means}
 &\nonumber
   \int_0^T\!\!\varPsi(\wt z)+\big\langle[\varPhi_0]_z'(\overline\varSigma_{\tau h},\underline{\overline z}_{\tau h}),\wt z{-}
    \DT z_{\tau h}\big\rangle_{\calZ_1^*\times\calZ_1}^{}\!
    +\big\langle\varPhi_1'\underline{\overline z}_{\tau h},\wt z\big\rangle_{\calZ^*\times\calZ}^{}
    \,\d t+\varPhi_1(z_0)
    \\[-.7em]&\hspace*{24em}
    \ge\varPhi_1(z_{\tau h}(T))+
    \!\int_0^T\!\!\varPsi(\DT z_{\tau h})\,\d t\,
    \label{weak-z-disc}
\end{align}\end{subequations}
\mman{for all $\wt z\in L^1(0,T;\calZ)$.}
This is completed by \eqref{suggestion-1+}. 

It is further important that the equations in
\eqref{suggestion-1+} and the first equation in \eqref{suggestion-2+}
are linear, so that the weak convergence is sufficient for the
limit passage there. In particular, we use \eqref{Eh->E} and the
Lebesgue dominated-convergence theorem.

As to the weak convergence of \eqref{suggestion-1+} integrated in time
towards \eqref{suggestion-1} integrated in time, i.e.\ towards $\varSigma
=\mathfrak CEu+G$ as used in Definition~\ref{def}, we need to prove that
\begin{align}
\int_0^T\!\!\big\langle\varSigma_{\tau h}-G_{\tau h},X\big\rangle_{\calS\times\calS^*}
-\big\langle u_{\tau h},E_h^*\mathfrak C^*X\big\rangle_{\calH\times\calH^*}\d t
\to\int_0^T\!\!\big\langle\varSigma{-}G,X\big\rangle_{\calS\times\calS^*}
-\big\langle u,E^*\mathfrak C^*X\big\rangle_{\calH\times\calH^*}\d t
\label{conv-S=CEu}\end{align}
for any $X\in
\mman{L^1(\mman{0,T;}\calS^*)}$.
By \eqref{Eh->E}, we have also $E_h^*S\to E^*S$ in $\calH$ for any
$S\in\calE^*$, in particular for $S=\mathfrak C^*X(t)$.
Thus certainly $E_h^*\mathfrak C^*X\to E^*\mathfrak C^*X$ in $L^1(\mman{0,T;}\calH)$
strongly. Using the weak* convergence $u_{\tau h}\to u$ in
$L^\infty(\mman{0,T;}\calH)$, we obtain \eqref{conv-S=CEu}.
Moreover, in the limit
$Eu=\mathfrak C^{-1}(\varSigma-G)\in L^\infty(\mman{0,T;}\calE)$ so that
$u\in L^\infty(\mman{0,T;}\calU)$.

For the limit passage in \eqref{weak-v-disc}, we also use
$\varPhi_\varSigma'(\overline\varSigma_{\tau h},\COLOR{\overline z_{\tau h}})\to
\varPhi_\varSigma'(\varSigma,z)$ weakly* in $L^\infty(\mman{0,T;}\calS^*)$
because $\varPhi_\varSigma'$ is continuous in the
(weak$\times$strong,weak)-mode, cf.~\eqref{ass-Phi-Sigma},
and because of the mentioned strong convergence of
$\overline z_{\tau h}\to z$.


Furthermore, we need to show the convergence 
$[\varPhi_0]_z'(\overline\varSigma_{\tau h},\underline{\overline z}_{\tau h})\to
[\varPhi_0]_z'(\varSigma,z)$.
For this, we use again the mentioned generalized Aubin-Lions theorem
to have the strong convergence $\underline{\overline z}_{\tau h}\to z$
in $L^{r}(\mman{0,T;}\calZ_1)$ for any $1\le r<+\infty$ and then
the continuity of $[\varPhi_0]_z'$ in the 
(weak$\times$strong,weak)-mode, cf.~the former option in \eqref{ass-Phi0-z}.
The limit passage of \eqref{weak-z-disc} towards
\eqref{weak-z} then uses also the weak lower semicontinuity of $\varPhi_1$ and
the weak convergence $z_{\tau h}(T)\to z(T)$ in $\calZ$; here for this
pointwise convergence in all time instants $t$ and in particular in $t=T$,
we also used that we have some information about $\DT z_{\tau h}$,
cf.\ \eqref{est3}.

So far, we have relied on the former options in (\ref{ass-Phi'}a,b) and the
Aubin-Lions compactness argument as far as the $z$-variable \cman{is concerned}.
If $\varPhi$ is quadratic (as e.g.\ in the examples in
Sects.\,\ref{sect-plast}--\ref{sect-poro} below), we can use the latter
options in (\ref{ass-Phi'}a,b) and simplify the above arguments, relying
merely on the weak convergence $\overline z_{\tau h}\to z$ and
$\overline{\underline z}_{\tau h}\to z$.
$\hfill\Box$

\begin{remark}[{\sl Alternative weak formulation}]
  \upshape
  Here, we used the weak formulation of \eqref{IVP2z+} containing the term
$\langle\varPhi_z'(\varSigma,z),\DT z\big\rangle$
which often does not have a good meaning since $\DT z$
may not be regular enough in some applications. This term is thus
eliminated
by substituting it, after integration over the time interval, by
$\varPhi(\varSigma(T),z(T))
-\int_0^T\langle\varPhi_\varSigma'(\varSigma,z),\DT\varSigma\big\rangle\,\d t
-\varPhi(\varSigma_0,z_0)$ or even rather by
$\varPhi(\varSigma(T),z(T))
-\int_0^T\langle\varPhi_\varSigma'(\varSigma,z),\mathfrak CEv\big\rangle\,\d t
-\varPhi(\varSigma_0,z_0)$. \cman{This} 
however, 
would bring even more difficulties
because we would need to prove
a strong convergence of
$\varPhi_\varSigma'(\varSigma,z)$, or of $\DT\varSigma$, or $\mathfrak CEv$ in
our explicit-discretization scheme, which seems not easy. 
\end{remark}

\begin{remark}[{\sl Nonquadratic $\varPhi(\varSigma,\cdot)$}]\label{rem-nonquadratic}
\upshape
Some applications use such $\varPhi(\varSigma,\cdot)$ which is not quadratic.
This is still consistent with the explicit leap-frog-type discretization if,
instead of $\varPhi_z'(\varSigma,z)$,
we consider an abstract difference quotient
$\varPhi_z^\circ(\varSigma,z,\widetilde z)$ with the properties
\begin{align}\label{quotient}
  \varPhi_z^\circ(\varSigma,z,z)=\varPhi_z'(\varSigma,z)
  \ \ \text{ and }\ \ \big\langle\varPhi_z^\circ(\varSigma,z,\widetilde z),
  z{-}\widetilde z\big\rangle
  =\varPhi(\varSigma,z)-\varPhi(\varSigma,\widetilde z)\,,
\end{align}
cf.\ \cite{RouPan17ECTD}. Then, instead of
$\varPhi_z'(\varSigma_{\tau h}^{\KK},\frac{z_{\tau h}^{k+1}{+}z_{\tau h}^{k}}{2})$ in
\eqref{suggestion-3}, \cman{we should} write
$\varPhi_z^\circ(\varSigma_{\tau h}^{\KK},z_{\tau h}^{k+1},z_{\tau h}^{k})$.
\end{remark}

\begin{remark}[{\sl State-dependent dissipation}]\label{rem-dissip}
\upshape
The generalization of $\varPsi$ dependent also on $z$ or even on $(\varSigma,z)$
is easy. Then $\partial\varPsi$ is to be replaced by the partial subdifferential
$\partial_{\DT z}\varPsi$ and \eqref{suggestion-3}
should use $\varPsi(\varSigma_{\tau h}^{\KK},z_{\tau h}^{k},\cdot)$
instead of $\varPsi(\cdot)$.
\end{remark}

\begin{remark}[{\sl Spatial numerical approximation}]\upshape
  From the coercivity of the stored energy $\varPhi$, we have
  $\varSigma_{\tau h}^k\in\calS$ for any $k=0,1,...$ and thus,
  from \eqref{suggestion-1},
  $E_hv_{\tau h}^k\in\calE$ so that $v_{\tau h}^k\in\calU$,
  although the limit $v$ cannot be assumed valued in $\calU$
  in general. Similarly, from \eqref{suggestion-2}, one can read that
  $E_h^*S_{\tau h}^k\in\calH$ although this cannot be expected in
  the limit in general. Anyhow, on the time-discrete level, one can use the
  FEM discretization similarly as in the linear elastodynamics where
  regularity can be employed, cf.\ \cite{BeJoTs02NFMF,BeJoTs01MFES,BeRoTs09CRFD,Tsog99} for a mixed finite\COLOR{-element} method and \cite{ChLaQi15SDGM} for the more recently developed staggered discontinuous Galerkin method for elastodynamics. 
\end{remark}

\mman{
  \begin{remark}[{\sl Other explicit-implicit schemes}]
    \upshape
  Combination of  explicit and implicit time discretization might not only be
  due to parabolic evolution of internal variables  but also due to
  geometrical reasons,  
            \cman{e.g. transmission through a thin layer,  
    that lead to a very restrictive} CFL condition, cf.\ \cite{BBDJ16ETCT}.
\end{remark}
}

  \section{Particular examples}\label{sec-exa}

We present \REPLACE{four}{three} examples from continuum mechanics of deformable
bodies at small strains of different characters to illustrate applicability
of the ansatz \eqref{ansatz} and the above discretization scheme.
Various combinations of these examples are possible, too, covering
thus a relatively wide variety of models.

We use a standard notation concerning function spaces.
Beside the Lebesgue $L^p$-spaces, we denote by $H^k(\varOmega;\R^n)$
the Sobolev space of functions whose distributional derivatives
are from $L^2(\varOmega;\R^{n\times d^k})$.

\subsection{Plasticity or creep}\label{sect-plast}

The simplest example with quadratic stored energy and local dissipation potential
is the model of plasticity or creep. The internal variable is then the plastic
strain $\pi$, valued in the set of symmetric trace-free matrices
$\R_{\rm dev}^{d\times d}=\{P\in\R^{d\times d};\ P^\top=P,\ {\rm tr}\,P=0\}$.
For simplicity, we consider only homogeneous Neumann or Dirichlet boundary
conditions, so that simply $E=e(u)$ and $\mathfrak C=\bbC$. 
The 
stored energy in terms of strain $e(u)$ is
\begin{align}
\mathscr{W}(u,\pi)=\int_\varOmega\frac12\mathbb C(e(u){-}\pi){:}
(e(u){-}\pi)\,\d x
\,,\label{W-creep}
\end{align}
which is actually a function of the elastic strain
$e_{\rm el}=e{-}\pi$. The additive decomposition $e(u)=e_{\rm el}{+}\pi$ is
referred to as Green-Naghdi's \cite{GreNag65GTEP} decomposition.
This energy leads to
\begin{align}
  \varPhi(\sigma,\pi)=\int_\varOmega\frac12\mathbb C^{-1}\sigma{:}\sigma
  -\sigma{:}\pi+\frac12\mathbb C\pi{:}\pi\,\d x
  \ \ \ \text{ with }\ \sigma=\mathbb C e(u)\,.
  \label{Phi-creep}\end{align} 
Let us note that $\varPhi_\sigma'(\sigma,\pi)=
\mathbb C^{-1}\sigma-\pi=e{-}\pi$, i.e.\ the elastic strain $e_{\rm el}$, and that
the proto-stress $\varSigma=\sigma$ is indeed different from the actual stress
$\sigma-\bbC\pi$.

The dissipation potential is standardly chosen as
\begin{align}
  \varPsi(\DT\pi)=\int_\varOmega\sigma_{_{\rm Y}}|\DT\pi|
  +\frac12\bbD\DT\pi{:}\DT\pi\,\d x
\end{align}
with $\sigma_{_{\rm Y}}\ge0$ a prescribed yield stress and $\bbD$ a positive
semidefinite viscosity tensor. The dissipation rate is then
$\varXi(\DT\pi)=\int_\varOmega\sigma_{_{\rm Y}}|\DT\pi|
  +\bbD\DT\pi{:}\DT\pi\,\d x$.
For $\bbD>0$ and $\sigma_{_{\rm Y}}=0$, we obtain mere creep model or, in
other words, the linear \emph{viscoelastic} model in the
\emph{Maxwell rheology}.
For both $\bbD>0$ and $\sigma_{_{\rm Y}}>0$, we obtain \emph{viscoplasticity}.
For $\bbD=0$ and $\sigma_{_{\rm Y}}>0$,
we would obtain the rate-independent (perfect) plasticity but
our Proposition~\ref{prop1} does not cover this case (i.e.\ $p=1$ is not
admitted).

The functional setting is $\calH=L^2(\varOmega;\R^d)$,
$\calE=\calS=\calZ=\calZ_1
=L^2(\varOmega;\R_{\rm sym}^{d\times d})$ where
$\R_{\rm sym}^{d\times d}$ denotes symmetric $(d{\times}d)$-matrices.
Thus $\calU:=\{v\In L^2(\varOmega;\R^d);\ e(v)\In L^2(\varOmega;\R^{d\times d})\}
=H^1(\varOmega;\R^d)$ by Korn's inequality.

A modification of the stored energy  models 
an {\it isotropic hardening}, enhancing \eqref{W-creep} as
\begin{align}
\mathscr{W}(u,\pi)=\int_\varOmega\frac12\mathbb C_1(e(u){-}\pi){:}
(e(u){-}\pi)+\frac12\mathbb C_2\pi{:}\pi\,\d x
\label{W-creep+}
\end{align}
so that the energy $\varPhi$ from \eqref{Phi-creep} is modified as
\begin{align}
\varPhi(\sigma,\pi)=\int_\varOmega\frac12\mathbb C_1^{-1}\sigma{:}\sigma
-\sigma{:}\pi+\frac12(\mathbb C_1{+}\mathbb C_2)\pi{:}\pi\,\d x\,.
\label{Phi-creep+}
\end{align}
In the pure creep variant $\sigma_{_{\rm Y}}=0$, this is actually the
{\it standard linear solid}
(in a so-called Zener form), considered together with the leap-frog
time discretization in \cite{BeEzJo04MFEA}. The isochoric constraint
${\rm tr}\,\pi=0$ can then be avoided, assuming that $\mathbb C_2$ is
positive definite.

All these models
lead to a flow rule which is
localized on each element when 
 an element-wise constant approximation of $\pi$ is used, 
 and \mman{no large system of algebraic equations need to be solved so that}
 the combination with the explicit discretization
of the other equations leads to a very fast computational procedure.

Another modification for {\it gradient plasticity} by adding terms
$\frac12\kappa|\nabla\pi|^2$ into the stored energy is easily
possible, too. This modification uses
$\calZ=H^1(\varOmega;\R_{\rm sym}^{d\times d})$ and
\eqref{z-ansatz} with
$\varPhi_1(z)=\int_\varOmega\frac12\kappa|\nabla\pi|^2$ and
makes, however, the flow rule nonlocal but at least \mman{one} can benefit
from that the usual space discretization
of the proto-stress $\sigma$ uses the continuous piecewise smooth elements
which allows for handling gradients $\nabla\pi$ if used consistently also for
$\pi$. 
For the quasistatic variant of this model, we refer to the classical
monographs \cite{HanRed99PMTN,Tema85MPP}, while the dynamical model with
$\bbD=0$ is e.g.\ in \cite[Sect.5.2]{MieRou15RIST}.

Noteworthy, all these models bear time regularity if the loading is smooth
and initial conditions regular enough, which can be advantageously reflected
in space FEM approximation, too. 

\mman{The Maxwell visco-elastodynamics was also studied by 
  J.-P.\,Groby \cite[Part~I, Sect.2]{Grob05MPOE} using a slightly modified time
  discretization scheme, namely the order of \eqref{suggestion-1}
and  \eqref{suggestion-2} was exchanged.}

\ccc
The CFL condition \eqref{CFLnew} here is actually the same as the standard
one \eqref{CFLold}. This is because the internal variable actually does not
influence the elasticity response and, likewise, the inertia is independent
of the internal variable, so the wave speed is not influenced either. 
The CFL  is thus of the form $\tau\le
h\sqrt{\varrho/|\lambda_{\rm max}(\bbC)|}$ where $\lambda_{\rm max}(\bbC)$ is the
maximal eigenvalue of $\bbC$. \eee
%

\subsection{Poroelasticity in isotropic materials}\label{sect-poro}

Another example with quadratic stored energy but less trivial dissipation
potential is a saturated {\it Darcy or Fick flow} of a diffusant in porous
media, e.g.\ water in porous elastic rocks or in concrete, or a solvent in
elastic polymers. The most simple model is the classical {\it Biot model}
\cite{Biot41GTTS}, capturing effects as swelling or seepage.
In \mman{a} one-component flow, the internal variable is then the
scalar-valued diffusant content (or concentration) denoted by $\zeta$.

 As in the previous \ccc Section~\ref{sect-plast} \eee, we consider only Neumann or
 Dirichlet boundary conditions, so that $E=e(u)$. Here we use the orthogonal
 decomposition $e={\rm sph}\,e+{\rm dev}\,e$
 with the spherical (volumetric) part ${\rm sph}\,e:=({\rm tr}\,e)\mathbb I/d$
 and the deviatoric part ${\rm dev}\,e$ and 
 confine ourselves to isotropic materials where the elastic-moduli tensor
 $\mathbb C_{ijkl}=K\delta_{ij}\delta_{kl}+G(\delta_{ik}\delta_{jl}
 +\delta_{il}\delta_{jk}-2\delta_{ij}\delta_{kl}/d)$ with
 $K$ the bulk modulus and $G$ the shear modulus (=\,the second Lam\'e
 constant), which is the standard notation hopefully without any
   confusion with the notation used in \eqref{IVP+}.
 Then the proto-stress $\varSigma=\sigma=\mathbb Ce=K{\rm sph}\,e+2G{\rm dev}\,e$.
 In particular, ${\rm sph}\,\sigma=K{\rm sph}\,e$ so that
 ${\rm tr}\,e=K^{-1}{\rm tr}\,\sigma$.
 
 Adopting the gradient theory for this internal variable $\zeta$, the 
 stored energy in terms of strain is considered
 \begin{align}\nonumber
   \mathscr{W}(u,\zeta)&=\int_\varOmega\frac12\mathbb Ce(u){:}e(u)+
   \frac12M(\beta{\rm tr}\,e(u){-}\zeta)^2
   +\frac12L(\zeta{-}\zeta_{\rm eq})^2+\frac\kappa2|\nabla\zeta|^2\,\d x
   \\&\nonumber
   =\int_\varOmega
   \frac12\Big(K+\frac{\beta^2}{d}M\Big)|{\rm sph}\,e(u)|^2+G|{\rm dev}\,e(u)|^2
   \\[-.3em]&\nonumber\qquad\qquad\qquad-\beta M\zeta{\rm tr}\,e(u)+\frac12M\zeta^2+\frac12L(\zeta{-}\zeta_{\rm eq})^2+\frac\kappa2|\nabla\zeta|^2\,\d x
 \label{W-poro}\end{align}
 which, in terms of the (here partial) stress $\sigma=\mathbb C e$,
 reads as
$\int_\varOmega\frac12(\frac1K+\frac{\beta^2}{dK^2}M)|{\rm sph}\,\sigma|^2
+\frac1G|{\rm dev}\,\sigma|^2-\beta\zeta\frac MK{\rm tr}\,\sigma+\frac12M\zeta^2
   +\frac12L(\zeta{-}\zeta_{\rm eq})^2\,\d x$.
Here $M>0$ and $\beta>0$ are so-called Biot modulus and coefficient,
respectively, $\kappa>0$ is a capillarity \REPLACE{constant}{coefficient},
and $\zeta_{\rm eq}$ is a given equilibrium content. 
 \mman{From \eqref{W-poro},} we arrive at the overall stored energy as:
\begin{align}\nonumber
  \varPhi(\sigma,\zeta)=
  \int_\varOmega&\frac12\Big(\frac1K+\frac{\beta^2}{dK^2}M\Big)|{\rm sph}\,\sigma|^2
  +\frac1G|{\rm dev}\,\sigma|^2-\beta\zeta\frac MK{\rm tr}\,\sigma\,\d x
  \\[-.3em]&\qquad\qquad\qquad\qquad+
  \!\!\!\!\!\lineunder{\int_\varOmega\frac12M\zeta^2
   +\frac12L(\zeta{-}\zeta_{\rm eq})^2
  +\frac\kappa2|\nabla\zeta|^2\,\d x\,,\!}{$=:\Phi_1(\zeta)$}\,.
\end{align}
Let us note that $\varPhi_\sigma'(\sigma,\zeta)=\bbC^{-1}\sigma+\frac{\beta M}{dK^2}(\beta{\rm sph}\,\sigma
{-}\zeta K\mathbb I)$, i.e.\ the elastic strain, and that
the proto-stress $\varSigma=\sigma$ indeed differs from an actual stress
by the spherical pressure part $\frac{\beta M}{dK}(\beta{\rm sph}\,\sigma
{-}\zeta K\mathbb I)$. 

The driving force for the diffusion is the
chemical potential $\mu=\varPhi_\zeta'(\sigma,\zeta)$, i.e.\ here
\begin{subequations}\label{CH}\begin{align}
  \mu=(M+L)\zeta-\beta \frac MK{\rm tr}\,\sigma
  -L\zeta_{\rm eq}-\kappa\Delta\zeta\,.
\end{align}
The diffusion equation is
\begin{align}
  \DT\zeta-{\rm div}(\mathbb M\nabla\mu)=0
\end{align}\end{subequations}
with $\mathbb M$ denoting the diffusivity tensor.
The system \eqref{CH} is called the {\it Cahn-Hilliard equation}, here combined
with elasticity so that the flow of the diffusant is driven both
by the gradient of concentration (Fick's law) and the gradient of
the mechanical pressure (Darcy's law). The dissipation potential in terms of
$\nabla\mu$, let us denote it by $R$ behind this system\mman{,} is
\begin{align}
R(\mu)=\int_\varOmega\frac12\mathbb M\nabla\mu{\cdot}\nabla\mu\,\d x\,.
\end{align}
For the analysis cf.\ e.g.\ \cite[Sect.\,7.6]{KruRou19MMCM}.

One would expect the dissipation potential as a function of the rate
of internal variables, as in \eqref{IVP2z+}. In fact,
the system \eqref{CH} turns into the form \eqref{IVP2z+} \cman{if} one takes
the dissipation potential $\varPsi=\varPsi(\DT\zeta)$ as
\begin{align}
\varPsi(\DT\zeta)=R^*(\DT\zeta)
\end{align}
with $R^*$ denoting the convex conjugate \cman{of} $R$. Now,
$\varPsi$ is nonlocal. 
The functional setting is as in the previous example but now
$\calZ=H^1(\varOmega)$ and $\calZ_1=H^1(\varOmega)^*$.
For a discretization of the
type \eqref{suggestion-3}, see \cite{Roub17ECTD}.

Often, the diffusivity is considered dependent on $\zeta$.
Or even one can think about $\mathbb M=\mathbb M(\sigma,\zeta)$. Then the
modification in Remark~\ref{rem-dissip} is to be applied. In particular,
$R(\sigma,\zeta,\mu)=\int_\varOmega\frac12\mathbb M(\sigma,\zeta)\nabla\mu
{\cdot}\nabla\mu\,\d x$ and $\varPsi(\sigma,\zeta,\DT\zeta)
=[R(\sigma,\zeta,\cdot)]^*(\DT\zeta)$.

  For this Biot model in the dynamical variant, the reader is also referred to
  the books \cite{AbChUl05PBC,Carc15WFRM,Chen16P,Stra08SWMP} or also
  \cite{KruRou19MMCM,MieRou15RIST}. In any case, the diffusion involves
  gradients and in the implicit discretization it leads to large systems of
  algebraic equations, which inevitably slows down the fast explicit
  discretization of the mechanical part itself.
  
  \ccc For this case also the CFL condition \eqref{CFLnew} is the same as the standard
one and leads to a restriction of the form $\tau \le C h/V_{\rm max}$ where $V_{\rm max}$ denotes the maximal speed with which waves propagate in the medium. 
We note that the pressure velocity which is the maximal speed of propagation in isotropic solids is enhanced in the Biot model. The stability analysis of the discrete scheme is 
quite technical and does not always lead to a practical CFL condition.
A.\,Ezziani in his Ph.D thesis \cite{EzzianiPhD} studied a discretization of
Biot's model similar to
\eqref{W-poro}
but the stability analysis of the discrete scheme
is very nontrivial, cf.\ 
\cite[Formula (7.4.11)]{EzzianiPhD}
and, as
he points out,
cannot be translated into a practical condition. 
  Therefore, in practice he proposes to use $\tau \le a_r h/V_{\rm pf}$
  where $V_{\rm pf}$ is the speed of the fast wave
  and $a_r$ is a constant
  depending on the order of the particular space discretization used.
  The attenuation caused by diffusion causes also some dispersion of
  wave velocities which stay however bounded from above by a high-frequency
  limit, cf.\ also \cite[Fig.\,5.2.1]{EzzianiPhD}, so the CFL condition
  expectedly holds uniformly like for the pure elastodynamics.
\eee
    
 \subsection{Damage}

The simplest examples of nonconvex stored energy are models of damage.
The most typical models use as an 
internal variable the scalar-valued bulk damage $\alpha$ having the
interpretation as a phenomenological 
volume fraction of microcracks or microvoids manifested macroscopically
as a certain weakening of the elastic response.
This concept was invented by L.M.\,Kachanov \cite{Kach58TRPD}
and Yu.N.\,Rabotnov \cite{Rabo69CPSM}.

Considering gradient theories, the stored energy in terms of the strain and
damage is here considered as
\begin{align}
\mathscr W(e,\alpha)
=\int_\varOmega\frac12\gamma(\alpha)\mathbb C e{:}e+\phi(\alpha)
+\frac\kappa2|\nabla\alpha|^2+\frac\varepsilon2\nabla(\mathbb C e)
{:}\nabla e\,\d x\,,
\label{W-dam}\end{align}
where $\phi(\cdot)$ is an energy of damage which gives rise to an activation
threshold for damage evolution and may also lead to healing (if allowed).
The last term is mainly to facilitate the mathematics towards convergence
and existence of a weak solution in such purely elastic materials
without involving any viscosity, cf.\
\cite[Sect.\,7.5.3]{KruRou19MMCM}. This regularization can also control
dispersion of elastic waves. \ccc More specifically, the 4th-order term
resulted in the momentum-equilibrium equation from the $\varepsilon$-term in 
\eqref{W-dam} causes an anomalous dispersion, i.e.\ waves with shorter
wavelength propagate faster than longer wavelength ones, cf.\ e.g.\
\cite[Remark 6.3.6]{KruRou19MMCM}. \eee
The $\nabla\alpha$-term also facilitates the analysis and controls
the internal length-scale of damage profiles.

Let us consider the ``generalized'' elasticity tensor $\mathfrak C=\bbC$
independent of $x$. As in the previous examples, $Eu=e(u)$ and $G=0$.
According \eqref{IVP2Sigma+}, the proto-stress $\varSigma=\mathfrak CEu+G$,
denoted by $\sigma$, now looks as $\mathbb C e=:\sigma$; in damage
  mechanics, the proto-stress is also called an effective stress with
  a specific mechanical interpretation, cf.\ \cite{Rabo69CPSM}.
In terms of $\sigma$, the stored energy is then
  \begin{align}
  \varPhi(\sigma,\alpha)=
  &\int_\varOmega\frac12\gamma(\alpha)\mathbb C^{-1}\sigma{:}\sigma
  +\frac\varepsilon2\nabla\mathbb C^{-1}\sigma{:}\nabla\sigma\,\d x
 +\!\!\!\!\lineunder{\int_\varOmega
  \phi(\alpha)+\frac\kappa2|\nabla\alpha|^2\,\d x\,.\!}{$=:\Phi_1(\alpha)$}
  \label{damage-stored}
  \end{align}
Then $\varPhi_\sigma'=\gamma(\alpha)\bbC^{-1}\sigma
    -{\rm div}(\varepsilon\nabla(\bbC^{-1}\sigma))$ and
    the true stress $S=\bbC^*\varPhi_\sigma'$ is then
    $\gamma(\alpha)\sigma-{\rm div}(\varepsilon\nabla\sigma)$ provided $\bbC$ is
    constant and symmetric.
 The damage driving force (energy) is $\varPhi_\alpha'(\sigma,\alpha)=
    \frac12\gamma'(\alpha)\mathbb C^{-1}\sigma{:}\sigma+\phi'(\alpha)
    -{\rm div}(\kappa\nabla\alpha)$. When $\gamma'(0)=0$ and $\phi'(0)\le0$,
    then always $\alpha\ge0$ also in the discrete scheme if $\alpha_0\ge0$.

    The other ingredient is the dissipation potential. To comply with the
    coercivity on $\calZ_1=L^2(\varOmega)$ with $p\ge2$ as needed in
    Proposition~\ref{prop1}, one can consider either
 \begin{align}\label{Psi-damage}  
   \varPsi(\DT\alpha)=\begin{cases}\int_\varOmega\varepsilon_1\DT\alpha^2\,\d x&
     \\\ \ +\infty&\end{cases}
     \text{ or }\ \ \ \begin{cases}\int_\varOmega\varepsilon_1\DT\alpha^2\,\d x&\text{if }\DT\alpha\le0\text{ a.e.\ on }\varOmega,
     \\\int_\varOmega\DT\alpha^2/\varepsilon_1\,\d x&\text{otherwise}\end{cases}
 \end{align}
 with some (presumably small) coefficient $\varepsilon_1>0$. The former option
 corresponds to a unidirectional (i.e.\ irreversible) damage not allowing any
 healing (as used in engineering) while the latter option allows for
 (presumably slow) healing as used in geophysical models on large time scales.
 
    Since $\sigma$ appears nonlinearly in $\varPhi_\alpha'(\sigma,\alpha)$,
the strong convergence $\overline\sigma_{\tau h}\to\sigma$ in $L^2(Q;\R^{d\times d})$
    is needed. For this, the strain-gradient term with $\varepsilon>0$ is needed
    and the Aubin-Lions compact embedding theorem is used. This gives
    the strong convergence even in the norm of
    $L^{1/\epsilon}(\mman{0,T;}L^{2d/(d-2)-\epsilon}(\varOmega;\R^{d\times d}))$ for
    arbitrarily small $\epsilon>0$ provided also $\DT\sigma_{\tau h}$ is bounded
    in some norm, which can be shown by using
    $\DT\sigma_{\tau h}=\bbC e(\overline v_{\tau h})$ and the Green formula
    \begin{align}\nonumber
      \big\|\DT\sigma_{\tau h}\|_{L^\infty(\mman{0,T;}H^{-1}(\varOmega;\R^{d\times d}))}^{}
      &=\sup_{\|\widetilde e\|_{L^1(\mman{0,T;}H_0^1(\varOmega;\R^{d\times d}))}\le1}
      \int_0^T\!\!\int_\varOmega\DT\sigma_{\tau h}{:}\widetilde e\,\d x\d t
      \\&\nonumber=\sup_{\|\widetilde e\|_{L^1(\mman{0,T;}H_0^1(\varOmega;\R^{d\times d}))}\le1}
     \int_0^T\!\!\int_\varOmega\bbC e(\overline v_{\tau h}){:}\widetilde e\,\d x\d t
      \\&\nonumber=\sup_{\|\widetilde e\|_{L^1(\mman{0,T;}H_0^1(\varOmega;\R^{d\times d}))}\le1}\!\!
-\int_0^T\!\!\int_\varOmega\overline v_{\tau h}{\cdot}{\rm div}(\bbC\widetilde e)\,\d x\d t
      \le C\|\overline v_{\tau h}\|_{L^\infty(\mman{0,T;}L^2(\varOmega;\R^d))}^{}
    \end{align}
    with $C$ depending on $|\bbC|$. Cf.\ also the abstract estimation
    \eqref{dual-est1}.

When $\gamma$ or $\phi$ are not quadratic but continuously
differentiable, one can use the abstract difference quotient \eqref{quotient}
defined, in the classical form, as
\begin{align}
  \varPhi_z^\circ(\varSigma,\alpha,\widetilde\alpha)=
  \begin{cases}\displaystyle{\frac12\frac{\gamma(\alpha){-}\gamma(\widetilde\alpha)}{\alpha{-}\widetilde\alpha}
    \mathbb C^{-1}\sigma{:}\sigma
    +\frac{\phi(\alpha){-}\phi(\widetilde\alpha)}{\alpha{-}\widetilde\alpha}
    -\kappa\Delta\frac{\alpha{+}\widetilde\alpha}2}
      &\text{where }\ \alpha\ne\widetilde\alpha\,.
      \\\frac12\gamma'(\alpha)\mathbb C^{-1}\sigma{:}\sigma+\phi'(\alpha)
        -\kappa\Delta\alpha&\text{where }\ \alpha=\widetilde\alpha\,.
    \end{cases}\label{nonquadratic-damaga}
\end{align}
Of course, rigorously, the $\Delta$-operator in \eqref{nonquadratic-damaga}
is to be understood in the weak form when using it in \eqref{suggestion-3}.

Due to the gradient $\kappa$-term in \eqref{damage-stored}, the implicit
incremental problem \eqref{suggestion-3} leads to an algebraic problem
with a full matrix, which may substantially slow down the otherwise fast
explicit scheme. Like in the previous model the capillarity, now this gradient
theory controls the length-scale of the damage profile and also serves as
a regularization to facilitate mathematical analysis. Sometimes, a nonlocal
``fractional'' gradient
can facilitate the analysis, too.
Then, some
wavelet equivalent norm can be considered to accelerate the calculations,
cf.\ also \cite{ArGrRo03MNSM}.
  As far as the stress-gradient term,
  it is important that the discretization of the proto-stress
  in the usual implementation of the leap-frog method is continuous piecewise
  smooth, so that $\nabla\sigma$ has a good sense in the discretization
  without need to use higher-order elements.
    Here we use that the latter relation in \eqref{suggestion-2} is to be
    understood in the weak form, namely
  $\int_\varOmega S_{\tau h}^{\KK}{:}\tilde E_h\,\d x
  =\langle\varPhi_\varSigma'(\varSigma_{\tau h}^{\KK},z_{\tau h}^{\KK}),
  \mathfrak C\tilde E_h\rangle$ for
  $\tilde E_h=\tilde e_h=e(\tilde u_h)$, which means
  $$
  \int_\varOmega S_{\tau h}^{\KK}{:}\tilde E_h\,\d x
  =\int_\varOmega\gamma(\alpha_{\tau h}^{k+1})\mathbb C^{-1}\sigma^{\KK}{:}
  \mathbb C\tilde e_h
  +\epsilon\nabla\mathbb C^{-1}\sigma^{\KK}{:}\nabla\mathbb C\tilde e_h\,\d x
  $$
  for any $\tilde e_h$ from the corresponding finite-dimensional subspace of
  $H^1(\varOmega;\R_{\rm sym}^{d\times d})$.
  Thus we indeed do not need higher-order elements, and also we do not need
  to specify explicitly homogeneous boundary conditions in this boundary-value problem.

  The functional setting is $\calH=L^2(\varOmega;\R^d)$, $\calE=\calS
  =H^1(\varOmega;\R_{\rm sym}^{d\times d})$,
  $\calZ=H^1(\varOmega)$, and $\calZ_0=\calZ_1=L^2(\varOmega)$.
  Then $\calU=H^2(\varOmega;\R^d)$, and $E=e(\cdot)$
  is understood as an operator $H^2(\varOmega;\R^d)\to
  H^1(\varOmega;\R_{\rm sym}^{d\times d})$, and $\mathfrak C^*\cong\bbC^\top\!=\bbC$
  is understood as an a operator from $H^1(\varOmega;\R_{\rm sym}^{d\times d})$
  to itself.


  \def\Frakg{{\mathfrak{g}}}
  A particular case of this model is a so-called {\it phase-field fracture},
  taking as \mman{a} basic choice
  \begin{align}\gamma(\alpha):={\varepsilon^2}/{\varepsilon_0^2}
        {+}\alpha^2,\ \ \ \
        \phi(\alpha):=\Frakg_{\rm c}{(1{-}\alpha)^2}/{\varepsilon},\ \
        \text{ and }\ \ \kappa:=\varepsilon \Frakg_{\rm c}
        \label{AT}\end{align}
  with $\Frakg_{\rm c}$ denoting the energy of fracture and with $\varepsilon$
  controlling a ``characteristic'' width of the {\it phase-field fracture}.
  The physical dimension of $\varepsilon_0$ as well as of 
$\varepsilon$ is m (meters) while the physical dimension of $\Frakg_{\rm c}$ 
is J/m$^2$. This is known as the so-called \emph{Ambrosio-Tortorelli 
  functional} \cite{AmbTor90AFDJ}. In the dynamical context, only various
implicit discretization schemes seems to be used so far, cf.\
\cite{BVSH12,HofMie12CPFM,RouVod19MMPF,SWKM14PFAD}.
There are a lot of improvements of this basic model, approximating 
a mode-sensitive fracture, or $\varepsilon$-insensitive models
(with $\varepsilon$ referring to \eqref{AT}), or ductile fracture\mman{,
  cf.\ \cite{Roub17ECTD}}. This last
variant combines this model with the plasticity as in Sect.~\ref{sect-plast}.

\ccc As mentioned above in this case we have anomalous dispersion, i.e.\ the high frequencies propagate faster, cf.\ e.g.\
\cite[Remark 6.3.6]{KruRou19MMCM}. The resulting CFL condition is a combination of the  
usual CFL  \eqref{CFLold} for the 2nd-order elastodynamic
model with the CFL condition for 4th-order plate as in
\cite{BeDeJo05ENMR}. More specifically,
the speed of elastic waves in such combined model is like 
$v\sim v_0\sqrt{1{+}\varepsilon/\lambda^2}$
with $v_0$ the speed in the elastodynamic case (i.e.\ $\varepsilon=0$)
and with $\lambda$ the wavelength, cf.\ \cite[Remark 6.3.6]{KruRou19MMCM} for a
one-dimensional analysis. For particular space discretisations,
implementable wavelengths $\lambda$
are bounded from below just by $h$. This yields to a CFL condition of the
type
\begin{align}
\tau\le C \frac h{\sqrt{1{+}\varepsilon/h^2}}.
\label{CFL-damage}\end{align}
Asymptotically, for
$h\to0$ we can see that $\tau$ is to be small as
$\mathscr{O}(\varepsilon^{-1/2}h^2)$.
For fixed $\varepsilon>0$, this is actually very restrictive like in the
explicit discretization
of the heat equation where it practically prevents from efficient
usage of explicit discretizations. However, here the role of $\varepsilon$ 
is primarily to facilitate rigorous existence of weak solutions of this
model and can be assumed to be small. Then the influence of
this 4th-order term and this restrictive asymptotics is presumably small, and
the usual CFL condition resulting from 
\eqref{CFL-damage} with $\varepsilon=0$ will dominate except
on very fine space discretizations.

Let us eventually remark that better asymptotics of the type 
$\tau\sim\mathscr{O}(\varepsilon^{-1/2}h^{1+\delta})$ for $h\to0$ can be obtained
by replacing the 4th-order term by a nonlocal term of the order
$2(1{+}\delta)$ for some $\delta>0$ small, which would
even allow
for more general dispersion \cite{Jira04NTCM} and simultaneously 
make the analytically desired regularization of the damage model,
cf.\ \cite[Remarks 6.3.7 and 7.5.29]{KruRou19MMCM}
\eee
\section*{Acknowledgments}



{This research has been partly supported from the grants 17-04301S
  (especially as far as the focus on the dissipative evolution of internal
  variables concerns)
and 19-04956S (especially as far as the focus on the dynamic and nonlinear
behaviour concerns)
of the Czech Sci.\ Foundation, and by the institutional support
RVO: 61388998 (\v CR). T.R.\ also acknowledges the hospitality of
FORTH in Heraklion, Crete.

\mman{The authors are deeply thankful \cman{to} Christos G. Panagiotopoulos for 
  fruitful discussions.
  \ccc Also the careful and critical reading of the original version
  and valuable suggestions by two anonymous referees are acknowledged. \eee
}
}

\baselineskip=11pt

\bibliographystyle{plain}
\bibliography{trcreta2}

\begin{thebibliography}{10}

\bibitem{AbChUl05PBC}
Y.N. Abousleiman, A.H.-D. Cheng, and F.-J. Ulm, editors.
\newblock {\em Poromechanics {III}: {B}iot Centennial (1905--2005)}, London,
  2005. Taylor \& Francis.

\bibitem{AmbTor90AFDJ}
L.~Ambrosio and V.M. Tortorelli.
\newblock Approximation of functional depending on jumps via by elliptic
  functionals via {$\Gamma$}-convergence.
\newblock {\em Comm. Pure Appl. Math.}, 43:999--1036, 1990.

\bibitem{ArGrRo03MNSM}
M.~Arndt, M.~Griebel, and T.Roub\'\i\v cek.
\newblock Modelling and numerical simulation of martensitic transformation in
  shape memory alloys.
\newblock {\em Continuum Mech. Thermodyn.}, 15:463--485, 2003.

\bibitem{ArnoldAwanou05}
D.N. Arnold and G.~Awanou.
\newblock Rectangular mixed finite elements for elasticity.
\newblock {\em Math. Models Methods Appl. Sci.}, 15:1417--1429, 2005.

\bibitem{ArnoldAwanouWinther08}
D.N. Arnold, G.~Awanou, and R.~Winther.
\newblock Finite elements for symmetric tensors in three dimensions.
\newblock {\em Math. Comput.}, 77:1229--1251, 2008.

\bibitem{ArnoldWinther02}
D.N. Arnold and R.~Winther.
\newblock Mixed finite elements for elasticity.
\newblock {\em Numerische Mathematik}, 92:401--419, 2002.

\bibitem{BeEzJo04MFEA}
E.~B\'ecache, A.~Ezziani, and P.~Joly.
\newblock A mixed finite element approach for viscoelastic wave propagation.
\newblock {\em Computational Geosciences}, 8:255--299, 2004.

\bibitem{BeDeJo05ENMR}
E.~B\'ecache, G.Derveaux, and P.~Joly.
\newblock An \ccc efficient numerical method for the resolution of the
  {K}irchhoff-{L}ove dynamic plate equation\eee.
\newblock {\em Numer. Meth. P. D. E.}, 21:323--348, 2005.

\bibitem{BeJoTs01MFES}
E.~B\'ecache, P.~Joly, and C.~Tsogka.
\newblock Fictitious domains, mixed finite elements and perfectly matched
  layers for {2D} elastic wave propagation.
\newblock {\em J. of Comput. Acoustics}, 9(3):1175--1202, 2001.

\bibitem{BeJoTs02NFMF}
E.~B\'ecache, P.~Joly, and C.~Tsogka.
\newblock A new family of mixed finite elements for the linear elastodynamic
  problem.
\newblock {\em SIAM J. Numer. Anal.}, 39:2109--2132, 2002.

\bibitem{BeRoTs09CRFD}
E.~B\'ecache, J.~Rodr\'{\i}guez, and C.~Tsogka.
\newblock Convergence \mman{results of the fictitious domain method for a mixed
  formulation of the wave equation with a {N}eumann boundary condition}.
\newblock {\em ESAIM: Math. Model. Numer. Anal.}, 43:377--398, 2009.

\bibitem{Biot41GTTS}
M.A. Biot.
\newblock General theory of three-dimensional consolidation.
\newblock {\em J. Appl. Phys.}, 12:155--164, 1941.

\bibitem{Bohl02PVFD}
T.~Bohlen.
\newblock Parallel {3-D} viscoelastic finite difference seismic modelling.
\newblock {\em Computers \& Geosciences}, 28:887--899, 2002.

\bibitem{BBDJ16ETCT}
M.~Bonnet, A.~Burel, M.~Durufl\'e, and P.~Joly.
\newblock Effective \mman{transmission conditions for thin-layer transmission
  problems in elastodynamics. {T}he case of a planar layer model}.
\newblock {\em ESAIM: Math. Model. Numer. Anal.}, 50:43--75, 2016.

\bibitem{BVSH12}
M.J. Borden, C.V. Verhoosel, M.A. Scott, T.J.R. Hughes, and C.M. Landis.
\newblock A phase-field description of dynamic brittle fracture.
\newblock {\em Comput. Meth. Appl. Mech. Engr.}, 217—-220:77--95, 2012.

\bibitem{Carc15WFRM}
J.M. Carcione.
\newblock {\em Wave Fields in Real Media, Wave Propagation in Anisotropic,
  Anelastic, Porous and Electromagnetic Media}.
\newblock Elsevier, Amsterdam, 2015.

\bibitem{Chen16P}
A.H.-D. Cheng.
\newblock {\em Poroelasticity}.
\newblock Springer, Switzerland, 2016.

\bibitem{ChLaQi15SDGM}
E.T. Chung, C.Y. Lam, and J.~Qian.
\newblock A staggered discontinuous {G}alerkin method for the simulation of
  seismic waves with surface topography.
\newblock {\em Geophysics}, 80:T119--T135, 2015.

\bibitem{CohPer17FEDG}
G.~Cohen and S.~Pernet.
\newblock {\em Finite Element and Discontinuous Galerkin Methods for Transient
  Wave Equations}.
\newblock Springer, Dordrecht, 2017.

\bibitem{CoFrLe28PDMP}
R.~Courant, K.~Friedrichs, and H.~Lewy.
\newblock \"{U}ber die partiellen {D}ifferenzengleichungen der mathematischen
  {P}hysik.
\newblock {\em Math. Annalen}, 100:32--74, 1928.

\bibitem{DelGli15AHOS}
S.~Delcourte and N.~Glinsky.
\newblock Analysis \mman{of a high-order space and time discontinuous
  {G}alerkin method for elastodynamic equations. {A}pplication to {3D} wave
  propagation}.
\newblock {\em ESAIM: Math. Model. Numer. Anal.}, 49:1085--1126, 2015.

\bibitem{EzzianiPhD}
A.~Ezziani.
\newblock {\em Mod\'elisation \ccc math\'ematique et num\'erique de la
  propagation d'ondes dans les milieux visco\'elastiques et
  poro\'elastiques\eee}.
\newblock PhD thesis, Univ.\ Paris IX Dauphine, 2005.

\bibitem{FaLeMa95MEIT}
C.~Farhat, M.~Lesoinne, and N.~Maman.
\newblock Mixed \mman{explicit/implicit time integration of coupled aeroelastic
  problems: Three‐field formulation, geometric conservation and distributed
  solution}.
\newblock {\em Intl. J. Numer. Meth. Fluids}, 21:807--835, 1995.

\bibitem{FePaFa01PACM}
C.A. Felippa, K.C. Park, and C.~Farhat.
\newblock Partitioned analysis of coupled mechanical systems.
\newblock {\em Comput. Methods Appl. Mech. Engrg.}, 190:3247--3270, 2001.

\bibitem{Grav96SSWP}
R.W. Graves.
\newblock Simulating seismic wave propagation in {3D} elastic media using
  staggered-grid finite differences.
\newblock {\em Bull. Seismological Soc. Amer.}, 86:1091--1106, 1996.

\bibitem{GreNag65GTEP}
A.E. Green and P.M. Naghdi.
\newblock A general theory of an elastic-plastic continuum.
\newblock {\em Arch. Rational Mech. Anal.}, 18:251--281, 1965.

\bibitem{Grob05MPOE}
J.-P. Groby.
\newblock {\em Mod\'elisation \mman{de la propagation des ondes \'elastiques
  g\'en\'er\'ees par un s\'eisme proche ou \'eloign\'e \`a l'int\'erieur d'une
  ville}}.
\newblock PhD thesis, Université de la M\'editerran\'ee - Aix-Marseille II,
  2005.

\bibitem{HanRed99PMTN}
W.~Han and B.D. Reddy.
\newblock {\em Plasticity}.
\newblock Springer, New York, 1999.

\bibitem{HofMie12CPFM}
M.~Hofacker and C.~Miehe.
\newblock Continuum phase field modeling of dynamic fracture: variational
  principles and staggered {FE} implementation.
\newblock {\em Intl. J. Fract.}, 178:113--129, 2012.

\bibitem{Jira04NTCM}
M.~Jir\'asek.
\newblock Nonlocal \ccc theories in continuum mechanics\eee.
\newblock {\em \ccc Acta Polytechnica\eee}, 44:16--34, 2004.

\bibitem{JolTso08FEMD}
P.~Joly and C.~Tsogka.
\newblock {\em Finite Element Methods with Discontinuous Displacement},
  chapter~11.
\newblock Chapman \& Hall/CRC, Boca Raton, FL, 2008.

\bibitem{Kach58TRPD}
L.M. Kachanov.
\newblock Time of rupture process under creep conditions.
\newblock {\em Izv. Akad. Nauk SSSR}, 8:26, 1958.

\bibitem{KPCG13GDAP}
R.~Kolman, J.~Ple\v{s}ek, J.~\v{C}erv, and D.~Gabriel.
\newblock Grid \mman{dispersion analysis of plane square biquadratic
  serendipity finite elements in transient elastodynamics}.
\newblock {\em Int. J. Numer. Meth. Engng.}, 96:1--28, 2013.

\bibitem{KPCOP16TSDS}
R.~Kolman, J.~Ple\v{s}ek, J.~\v{C}erv, M.~Okrouhl{\'\i}k, and
  P.~Pa\v{r}{\'\i}k.
\newblock Temporal-spatial \mman{dispersion and stability analysis of finite
  element method in explicit elastodynamics}.
\newblock {\em Intl. J. Numer. Meth. Engr.}, 106:113--128, 2016.

\bibitem{KruRou19MMCM}
M.~Kru{\v{z}}{\'\i}k and T.~Roub{\'\i}{\v{c}}ek.
\newblock {\em Mathematical Methods in Continuum Mechanics of Solids}.
\newblock Springer, Switzeland, 2019.

\bibitem{Maug15SIVS}
G.A. Maugin.
\newblock The saga of internal variables of state in continuum thermo-mechanics
  (1893-2013).
\newblock {\em Mechanics Research Communications}, 69:79--86, 2015.

\bibitem{MieRou15RIST}
A.~Mielke and T.~Roub{\'\i}{\v{c}}ek.
\newblock {\em Rate-Independent Systems -- Theory and Application}.
\newblock Springer, New York, 2015.

\bibitem{Rabo69CPSM}
Yu.N. Rabotnov.
\newblock {\em Creep Problems in Structural Members}.
\newblock North-Holland, Amsterdam, 1969.

\bibitem{Roub13NPDE}
T.~Roub{\'\i}{\v{c}}ek.
\newblock {\em Nonlinear Partial Differential Equations with Applications}.
\newblock Birkh\"auser, Basel, 2nd edition, 2013.

\bibitem{Roub17ECTD}
T.~Roub{\'\i}{\v{c}}ek.
\newblock An energy-conserving time-discretisation scheme for poroelastic media
  with phase-field fracture emitting waves and heat.
\newblock {\em Disc. Cont. Dynam. Syst. S}, 10:867--893, 2017.

\bibitem{RKVPZ??DACM}
T.~Roub{\'\i}{\v{c}}ek, M.~Kru\v{z}{\'\i}k, V.~Manti{\v{c}}, C.G.
  Panagiotopoulos, R.~Vodi\v{c}ka, and J.~Zeman.
\newblock Delamination and adhesive contacts, their mathematical modeling and
  numerical treatment.
\newblock In V.Manti\v{c}, editor, {\em Math. Methods and Models in
  Composites}, chapter~11. Imperial College Press, 2nd edition.

\bibitem{RouPan17ECTD}
T.~Roub{\'\i}{\v{c}}ek and C.G. Panagiotopoulos.
\newblock Energy-conserving time-discretisation of abstract dynamical problems
  with applications in continuum mechanics of solids.
\newblock {\em Numer. Funct. Anal. Optim.}, 38:1143--1172, 2017.

\bibitem{RoPaTs19ETDE}
T.~Roub{\'\i}{\v{c}}ek, C.G. Panagiotopoulos, and C.~Tsogka.
\newblock Explicit \mman{time-discretisation of elastodynamics with some
  inelastic processes at small strains}.
\newblock Preprint arXiv no.1903.11654, 2019.

\bibitem{RouVod19MMPF}
T.~Roub{\'\i}{\v{c}}ek and R.~Vodi\v{c}ka.
\newblock A monolithic model for phase-field fracture and waves in solid-fluid
  media towards earthquakes.
\newblock {\em Intl. J. Fracture}, 219, 2019.

\bibitem{Scar04ETNP}
G.~Scarella.
\newblock {\em Etude th\'eorique et num\'erique de la propagation d'ondes en
  pr\'esence de contact unilat\'eral dans un milieu fissur\'e}.
\newblock PhD thesis, Univ. Paris Dauphine, 2004.

\bibitem{SWKM14PFAD}
A.~Schl\"uter, A.~Willenb\"ucher, C.~Kuhn, and R.~M\"uller.
\newblock Phase field approximation of dynamic brittle fracture.
\newblock {\em Comput Mech}, 54:1141--1161, 2014.

\bibitem{SePaPa18SEIF}
S.~Seifi, K.C. Park, and H.S. Park.
\newblock A staggered explicit-implicit finite element formulation
  forelectroactive polymers.
\newblock {\em Comput. Methods Appl. Mech. Engrg.}, 337:150--164, 2018.

\bibitem{Stra08SWMP}
B.~Straughan.
\newblock {\em Stability and Wave Motion in Porous Media}.
\newblock Springer, New York, 2008.

\bibitem{Tema85MPP}
R.~Temam.
\newblock {\em Mathematical Problems in Plasticity}.
\newblock Gauthier-Villars, Paris, 1985.
\newblock (French original in 1983).

\bibitem{Tsog99}
C.~Tsogka.
\newblock {\em Modelisation math\'ematique et num\'erique de la propagation des
  ondes \'elastiques tridimensionnelles dans des milieux fissur\'es}.
\newblock PhD thesis, Univ. Paris IX Dauphine, 1999.

\bibitem{Viri84SHWP}
J.~Virieux.
\newblock {SH}-wave propagation in heterogeneous media: {Velocity}-stress
  {finite}-difference method.
\newblock {\em Geophysics}, 49:1933--1957, 1984.

\end{thebibliography}

\end{document}